\documentclass[12pt]{article}
\usepackage{theorem} 

\usepackage{amsfonts}
\usepackage{latexsym}
\usepackage{amssymb}
\usepackage{amsbsy}

\newtheorem{thm}{Theorem}

\newtheorem{prop}[thm]{Proposition} 

\newtheorem{lemma}[thm]{Lemma} 

\newtheorem{cor}[thm]{Corollary}

\newcommand{\bprf}[1][Proof:]{\begin{list}{} 			{\setlength{\leftmargin}{1em} 			\setlength{\rightmargin}{0em}}                         \item {\bf \hspace{-1em}  #1 \ \ }} 

\begin{document}

%\newlength{\oldbaselineskip}
%\setlength{\oldbaselineskip}{\baselineskip}
%\setlength{\baselineskip}{1.5\oldbaselineskip}

\title{Multiplicative Cellular Automata on Nilpotent Groups: 
Structure, Entropy, and Asymptotics}

\author{Marcus Pivato\\ {\em Department of Mathematics, Trent University}\\
 {\em Peterborough, Ontario, Canada}\\
 email: {\tt mpivato@trentu.ca} \\ Phone: (705) 748 1011 x1293  \ \ \ \  Fax: (705) 748 1630}

\maketitle

\begin{abstract}
  If ${\mathbb{M}}={\mathbb{Z}}^D$, and ${\mathcal{ B}}$ is a finite (nonabelian) group,
then ${\mathcal{ B}}^{\mathbb{M}}$ is a compact group; \ a {\em multiplicative
cellular automaton} (MCA) is a continuous transformation
${\mathfrak{ G}}:{\mathcal{ B}}^{\mathbb{M}}\!\longrightarrow{\mathcal{ B}}^{\mathbb{M}}\!$ which commutes with all shift maps, and where
nearby coordinates are combined using the multiplication operation of
${\mathcal{ B}}$.

  We characterize when MCA are group endomorphisms of ${\mathcal{ B}}^{\mathbb{M}}$, and
show that MCA on ${\mathcal{ B}}^{\mathbb{M}}$ inherit a natural structure theory from
the structure of ${\mathcal{ B}}$.  We apply this structure theory to compute
the measurable entropy of MCA, and to study convergence of initial
measures to Haar measure.

\end{abstract}

\section{Introduction}

  If ${\mathcal{ B}}$ is a finite set, and ${\mathbb{M}}$ is some indexing set, then the {\bf
configuration space} ${\mathcal{ B}}^{\mathbb{M}}$ is the set of all ${\mathbb{M}}$-indexed
sequences of elements on ${\mathcal{ B}}$.  If ${\mathcal{ B}}$ is discretely topologised,
then the Tychonoff product topology on ${\mathcal{ B}}^{\mathbb{M}}$ is compact, totally
disconnected, and metrizable.  If ${\mathbb{M}}$ is an abelian
monoid (e.g. ${\mathbb{M}} = {\mathbb{Z}}^D$, ${\mathbb{N}}^E$, or ${\mathbb{Z}}^D \times {\mathbb{N}}^E$), then the action of
${\mathbb{M}}$ on itself by translation induces a natural {\bf shift action} of
${\mathbb{M}}$ on configuration space: \ for all ${\mathsf{ v}} \in {\mathbb{M}}$, and ${\mathbf{ b}} =
{\left[b_{\mathsf{ m}}  |_{{\mathsf{ m}}\in{\mathbb{M}}}^{} \right]}
\in {\mathcal{ B}}^{\mathbb{M}}$, define $ {{{\boldsymbol{\sigma}}}^{{\mathsf{ v}}}} [{\mathbf{ b}}] \ = \ {\left[b'_{\mathsf{ m}}  |_{{\mathsf{ m}}\in{\mathbb{M}}}^{} \right]}$
 where, $\forall {\mathsf{ m}}, \ \ \ b'_{\mathsf{ m}} = b_{{\mathsf{ v}}+{\mathsf{ m}}}$.

 A {\bf cellular automaton} (CA) is a continuous map
${\mathfrak{ G}}:{\mathcal{ B}}^{\mathbb{M}}\!\longrightarrow{\mathcal{ B}}^{\mathbb{M}}\!$ which commutes with all shifts: \ for any
${\mathsf{ m}}\in{\mathbb{M}}$,\ \ \ $ {{{\boldsymbol{\sigma}}}^{{\mathsf{ m}}}} \circ
{\mathfrak{ G}} \ = \ {\mathfrak{ G}} \circ  {{{\boldsymbol{\sigma}}}^{{\mathsf{ m}}}} $.  A result of Curtis, Hedlund, and
Lyndon \cite{HedlundCA}
says any CA is determined by a {\bf local map} ${\mathfrak{ g}}:{\mathcal{ B}}^{\mathbb{V}} {{\longrightarrow}}
{\mathcal{ B}}$ (where ${\mathbb{V}}\subset {\mathbb{M}}$ is some finite subset), so that, for all
${\mathsf{ m}}\in{\mathbb{M}}$, if we define ${\mathsf{ m}}+{\mathbb{V}} = {\left\{ {\mathsf{ m}}+{\mathsf{ v}} \; ; \; {\mathsf{ v}}\in{\mathbb{V}} \right\} }$, and for all ${\mathbf{ b}}
\in
{\mathcal{ B}}^{\mathbb{M}}$, if we define ${\mathbf{ b}}\raisebox{-0.3em}{$\left|_{({\mathsf{ m}}+{\mathbb{V}})}\right.$}$ to be the restriction of ${\mathbf{ b}}$
to an element of ${\mathcal{ B}}^{({\mathsf{ m}}+{\mathbb{V}})}$, then ${\mathfrak{ G}}({\mathbf{ b}})_{\mathsf{ m}} \ = \
{\mathfrak{ g}}\left({\mathbf{ b}}\raisebox{-0.3em}{$\left|_{{\mathsf{ m}}+{\mathbb{V}}}\right.$}\right)$.
 
  If $({\mathcal{ B}},\cdot)$ is a finite multiplicative group,  let
${{{\mathbf{ E}}{\mathbf{ n}}{\mathbf{ d}}_{} \left[{\mathcal{ B}}\right]}}$ be the set of group endomorphisms of ${\mathcal{ B}}$.  A
{\bf multiplicative cellular automaton} (MCA) is a CA whose local map
is a product of affine endomorphisms of separate coordinates.  To be
precise, let $ v:{\left[ 1...I \right]}{{\longrightarrow}}{\mathbb{V}}$ be a (possibly
noninjective) map, let ${\mathfrak{ g}}_1,{\mathfrak{ g}}_2,\ldots,{\mathfrak{ g}}_I\in{{{\mathbf{ E}}{\mathbf{ n}}{\mathbf{ d}}_{} \left[{\mathcal{ B}}\right]}}$, and
let $ g_0, g_1,\ldots, g_I\in{\mathcal{ B}}$ be constants.  If ${\mathbf{ b}} =
{\left[b_{\mathsf{ v}}  |_{{\mathsf{ v}}\in{\mathbb{V}}}^{} \right]} \in {\mathbb{G}}^{\mathbb{V}}$, then the local map ${\mathfrak{ g}}$ has
the form:
\begin{equation}
\label{local.map.eqn}
  {\mathfrak{ g}}({\mathbf{ b}}) \ = \  g_0 \cdot {\mathfrak{ g}}_{1}\left(b_{ v[1]}\right) \cdot  g_1 \cdot
	{\mathfrak{ g}}_{2}\left(b_{ v[2]}\right) \cdot   g_2 \cdot
	\ldots \cdot   g_{I-1} \cdot {\mathfrak{ g}}_{I}\left(b_{ v[I]}\right) \cdot   g_I
\end{equation}
  The {\bf ordering function} $ v$ imposes an order on this product, which
is necessary if ${\mathcal{ B}}$ is nonabelian.  The endomorphisms
${\left[{\mathfrak{ g}}_i  |_{i=1}^{I} \right]}$ are called the {\bf coefficients} of ${\mathfrak{ G}}$.
We can rewrite equation (\ref{local.map.eqn}) as
\begin{equation}
\label{local.map.eqn.2}
 {\mathfrak{ g}}({\mathbf{ b}}) \ = \  g \cdot
  {\mathfrak{ g}}'_{1}\left(b_{ v[1]}\right)\cdot {\mathfrak{ g}}'_{2}\left(b_{ v[2]}\right) \cdot  
	\ldots \cdot {\mathfrak{ g}}'_{I}\left(b_{ v[I]}\right)
\ \ = \ \  g\cdot \prod_{i=1}^{I} {\mathfrak{ g}}'_{i}\left(b_{ v[i]}\right),
\end{equation}
where  $ g \ = \  g_0 \cdot  g_1 \cdot \ldots \cdot  g_I$, and,
for each $i\in{\left[ 1...I \right]}$, \ 
${\mathfrak{ g}}'_{i}(b) \ = \ \left( g_{I}  g_{I-1}\ldots  g_{i}\right)^{-1} \cdot {\mathfrak{ g}}_{i}(b)
\cdot \left( g_{I}  g_{I-1}\ldots  g_{i}\right)$
is an endomorphism.  The product ``$\prod_{i=1}^I$'' inherits the
obvious order from ${\left[ 1...I \right]}$.  We assume MCAs are written in the
form (\ref{local.map.eqn.2}), and call $ g$ the {\bf bias}.  If the
bias is trivial (${\mathfrak{ G}}$ is ``unbiased''), then ${\mathfrak{ g}}({\mathbf{ b}})$ is just a
product of endomorphic images of the components $\{b_{\mathsf{ v}}\}_{{\mathsf{ v}}\in{\mathbb{V}}}$.

  ${\mathcal{ B}}^{\mathbb{M}}$ is a
compact group under componentwise multiplication; \ an {\bf
endomorphic cellular automaton} (ECA) is a topological group endomorphism
${\mathfrak{ G}}:{\mathcal{ B}}^{\mathbb{M}}\!\longrightarrow{\mathcal{ B}}^{\mathbb{M}}\!$ which commutes with all shift maps.  If ${\mathcal{ B}}$ is
abelian, then all unbiased MCA are ECA, and vice versa;\ when ${\mathcal{ B}}$ is
nonabelian, however, the ECA form only a small subclass of
MCA (see \S\ref{S:ECA}).

  \medskip         \refstepcounter{thm}                     \begin{list}{} 			{\setlength{\leftmargin}{1em} 			\setlength{\rightmargin}{0em}}                     \item {\bf Example \thethm:} 
Consider the following local maps:
\setcounter{enumi}{1}
\begin{list}{(\alph{enumii})}{\usecounter{enumii}}
  \item \label{x0} Let ${\mathbb{M}}={\mathbb{Z}}$, \ ${\mathbb{V}}=\{0,1\}$, and let
${\mathfrak{ g}}\left(b_0,b_1\right) \ = \ b_0\cdot b_1$.
  \item \label{x1} ${\mathbb{M}}={\mathbb{N}}$, \ ${\mathbb{V}}={\left[ 0...2 \right]}$; \ 
${\mathfrak{ g}}\left(b_0,b_1,b_2\right) \ = \ b_0\cdot b_1\cdot b_2$.
  \item  \label{x2} ${\mathbb{M}}={\mathbb{N}}$, \ ${\mathbb{V}}={\left[ 0...2 \right]}$; \ 
 ${\mathfrak{ g}}\left(b_0,b_1,b_2\right) \ = \ b_2^4\cdot b_1^3\cdot b_0$.
  \item \label{x8} ${\mathbb{M}}={\mathbb{Z}}^2$, \  ${\mathbb{V}}={\left[ -1...1 \right]}^2$; \
${\mathfrak{ g}}({\mathbf{ b}}) \ = \ \left( g \cdot b_{(-1,0)}\cdot  g^{-1}\right)
 \cdot b_{(0,-1)}\cdot b_{(0,0)} \cdot  b_{(0,1)}\cdot h \cdot  b_{(1,0)}^{-1}$, where $ g,h\in{\mathcal{ B}}$ are constants.

  \item  \label{x7} Suppose ${\mathcal{ B}} = {\mathbb{G}}{\mathbb{L}}^n({\mathbb{F}})$ is the group of invertible $n \times n$ matrices over
a finite field ${\mathbb{F}}$ and let ${\mathfrak{ g}}\left( {\mathbf{ B}}_{-1}, {\mathbf{ B}}_0, {\mathbf{ B}}_1 \right) \ = \ 
\det[{\mathbf{ B}}_{-1}] \cdot \det[{\mathbf{ B}}_{1}]^2 \cdot {\mathbf{ B}}_0.$
\end{list}

  Example (\ref{x0}) is the  {\bf nearest-neighbour multiplication} CA
 \cite{MooreQuasi,MooreNLCA}.
  Examples (\ref{x0}-\ref{x2}) are unbiased, and all 
coefficients are the identity map on ${\mathcal{ B}}$.

In (\ref{x2}),\ $ v:{\left[ 1...8 \right]}{{\longrightarrow}}
{\left[ 0...2 \right]}$ is defined: $ v[1]= v[2]= v[3]= v[4]=2$, $ v[5]= v[5]= v[7]=1$, and $ v[8]=0$;\  by repeating indices in this way, we can obtain any
exponents we want. 

In (\ref{x8}), suppose ${{\sf card}\left[{\mathcal{ B}}\right]}= B$; then
$ v:{\left[ 1\ldots 4\!+\! B \right]}{{\longrightarrow}}{\mathbb{V}}$ is defined:
$ v[1]=(-1,0)$, $ v[2]=(0,-1)$, $ v[3]=(0,0)$, $ v[4]=(0,1)$, and $ v[n]=(1,0)$ for $n\in{\left[ 5\ldots 4\!+\! B \right]}$; in this way, we
obtain the exponent $b_{(1,0)}^{ B-1} = b_{(1,0)}^{-1}$.
All 
coefficients are the identity map, except for 
${\mathfrak{ g}}_1(b) =  g\cdot b\cdot  g^{-1}$, which is the
endomorphism of conjugation-by-$ g$.

In (\ref{x7}), let ${\mathbf{ I}}\in{\mathbb{G}}{\mathbb{L}}^n({\mathbb{F}})$ be the
identity matrix.  Then ${\mathfrak{ g}}_1({\mathbf{ B}}) = \det[{\mathbf{ B}}]\cdot{\mathbf{ I}}$
and ${\mathfrak{ g}}_2({\mathbf{ B}}) = \det[{\mathbf{ B}}]^2\cdot{\mathbf{ I}}$ are endomorphisms of ${\mathbb{G}}{\mathbb{L}}^n({\mathbb{F}})$,
and ${\mathfrak{ g}}_3={\mathbf{ Id}_{{}}}$.  In fact, (\ref{x7}) is an ECA.
  	\hrulefill\end{list}   			  

\medskip

   When ${\mathcal{ B}}$ is an additive abelian group
(e.g. ${\mathcal{ B}}=({{\mathbb{Z}}_{/p}},+)$), unbiased MCA are called {\bf linear} CA
(or {\bf affine} CA, when biased).  Classical modular arithmetic has
been applied to study the entropy
\cite{DamicoManziniMargara}, and computational complexity
\cite{MooreQuasi,MooreNLCA} of linear CA, while techniques of harmonic
analysis yield convergence of initial probability measures on
${\mathcal{ B}}^{\mathbb{M}}$ to the uniformly distributed, or {\bf Haar} measure under
iteration by affine CA
\cite{MaassMartinezII,Lind,PivatoYassawi1,PivatoYassawi2,FerMaassMartNey}.
However, the case when ${\mathcal{ B}}$ is nonabelian is poorly understood;
\ ``abelian'' techniques usually fail to apply.

\medskip

   In \S\ref{S:ECA}, we give necessary and sufficient conditions for
an MCA to be endomorphic.  In
\S\ref{S:structure}, we use the structure theory of the group ${\mathcal{ B}}$ to
develop a corresponding structure theory for MCA over ${\mathcal{ B}}$.  We apply
this structure theory in \S\ref{S:entropy}, to compute the measurable
entropy of MCA, and in \S\ref{S:measure} to establish sufficient
conditions for convergence of initial measures to Haar measure under
iteration of MCA.  The major results are
Theorems \ref{MCA.structure}, \ref{entropy.skewprod}, and \ref{nilpotent.to.haar}.

\section{Endomorphic Cellular Automata
\label{S:ECA}}

  Suppose ${\mathfrak{ G}}:{\mathcal{ B}}^{\mathbb{M}}\!\longrightarrow{\mathcal{ B}}^{\mathbb{M}}\!$ is an ECA.
Since ${\mathfrak{ G}}\in{{{\mathbf{ E}}{\mathbf{ n}}{\mathbf{ d}}_{} \left[{\mathcal{ B}}^{\mathbb{M}}\right]}}$, the local map ${\mathfrak{ g}}$ must
be a group homomorphism from the product group ${\mathcal{ B}}^{\mathbb{V}}$ into ${\mathcal{ B}}$.
This constrains the coefficients
$\{{\mathfrak{ g}}_{\mathsf{ v}}\}_{{\mathsf{ v}}\in{\mathbb{V}}}$ and their interactions.

\begin{lemma}\label{ECA.local.map}  
  Let ${\mathfrak{ g}}:{\mathcal{ B}}^{\mathbb{V}} {{\longrightarrow}} {\mathcal{ B}}$ be a group homomorphism.
Then there are endomorphisms
${\mathfrak{ g}}_{\mathsf{ v}}\in{{{\mathbf{ E}}{\mathbf{ n}}{\mathbf{ d}}_{} \left[{\mathcal{ B}}\right]}}$ for all ${\mathsf{ v}}\in{\mathbb{V}}$ so that, for any \
${\mathbf{ b}} \ = \ {\left[b_{\mathsf{ v}}  |_{{\mathsf{ v}}\in{\mathbb{V}}}^{} \right]} \in {\mathcal{ B}}^{\mathbb{V}}$, \ \ 
$\displaystyle {\mathfrak{ g}}({\mathbf{ b}}) \ = \ \prod_{{\mathsf{ v}}\in{\mathbb{V}}} {\mathfrak{ g}}_{\mathsf{ v}}(b_{\mathsf{ v}})$, 
where this product is commutative. \end{lemma}
\bprf
  For each ${\mathsf{ v}}\in{\mathbb{V}}$, let ${\mathfrak{ i}}_{\mathsf{ v}}:{\mathcal{ B}}{{\longrightarrow}}{\mathcal{ B}}^{\mathbb{V}}$ be the 
embedding into the ${\mathsf{ v}}$th coordinate: \ for any $b\in{\mathcal{ B}}$,
\  $\left({\mathfrak{ i}}_{\mathsf{ v}}(b)\right)_{\mathsf{ v}} = b$, and
$\left({\mathfrak{ i}}_{\mathsf{ v}}(b)\right)_{\mathsf{ w}} = e$ for all ${\mathsf{ w}}\neq {\mathsf{ v}}$ in ${\mathbb{V}}$, where $e\in{\mathcal{ B}}$
is the identity element.  Then define ${\mathfrak{ g}}_{\mathsf{ v}} = {\mathfrak{ g}}\circ{\mathfrak{ i}}_{\mathsf{ v}}$.  If
${\mathbf{ b}} = {\left[b_{\mathsf{ v}}  |_{{\mathsf{ v}}\in{\mathbb{V}}}^{} \right]}$, then clearly,
$\displaystyle {\mathbf{ b}} \ = \ \prod_{{\mathsf{ v}}\in{\mathbb{V}}} {\mathfrak{ i}}_{\mathsf{ v}}(b_{\mathsf{ v}})$,  \ 
where the factors all commute, and thus,
$\displaystyle {\mathfrak{ g}}({\mathbf{ b}}) 
\ = \ {\mathfrak{ g}}\left(\prod_{{\mathsf{ v}}\in{\mathbb{V}}} {\mathfrak{ i}}_{\mathsf{ v}}(b_{\mathsf{ v}}) \right) 
\ = \ \prod_{{\mathsf{ v}}\in{\mathbb{V}}} {\mathfrak{ g}}\left({\mathfrak{ i}}_{\mathsf{ v}}(b_{\mathsf{ v}}) \right)
\ = \ \prod_{{\mathsf{ v}}\in{\mathbb{V}}} {\mathfrak{ g}}_{\mathsf{ v}}(b_{\mathsf{ v}})$,
where, again, the factors all commute.
{\tt \hrulefill $\Box$ }\end{list}  \medskip  

  We say two endomorphisms ${\mathfrak{ g}}_{\mathsf{ w}}$ and ${\mathfrak{ g}}_{\mathsf{ v}}$ have {\bf commuting
images} if, for any $b_{\mathsf{ w}},b_{\mathsf{ v}} \in {\mathcal{ B}}$, \ \ ${\mathfrak{ g}}_{\mathsf{ v}}(b_{\mathsf{ v}})\cdot {\mathfrak{ g}}_{\mathsf{ w}}(b_{\mathsf{ w}})
\ = \ {\mathfrak{ g}}_{\mathsf{ w}}(b_{\mathsf{ w}})\cdot {\mathfrak{ g}}_{\mathsf{ v}}(b_{\mathsf{ v}})$.  Thus, the coefficients of any ECA
${\mathfrak{ G}}$ must all have commuting images; \ this  restricts
the structure of ${\mathfrak{ G}}$, and the more noncommutative ${\mathcal{ B}}$ itself
is, the more severe the restriction becomes.  The noncommutativity of
${\mathcal{ B}}$ is measured by two subgroups:  \ the {\bf centre}, $Z({\mathcal{ B}}) =
{\left\{ z\in{\mathcal{ B}} \; ; \; \forall b\in{\mathcal{ B}}, \ b\cdot z = z\cdot b \right\} }$, and the {\bf commutator
subgroup}, $[{\mathcal{ B}},{\mathcal{ B}}] =
 {\left\langle b\cdot h\cdot b^{-1}\cdot h^{-1}; \ b,h\in{\mathcal{ B}} \right\rangle }$.
If $\phi:{\mathcal{ B}}{{\longrightarrow}}{\mathcal{ A}}$ is any homomorphism from ${\mathcal{ B}}$ into an abelian group
${\mathcal{ A}}$, then $[{\mathcal{ B}},{\mathcal{ B}}]\subset\ker[\phi]$.
\begin{cor}\label{limitations.of.ECA.coefficients}    Contining with the previous notation,
\begin{enumerate}
\item If \ $\exists{\mathsf{ v}}\in{\mathbb{V}}$ so that ${\mathfrak{ g}}_{\mathsf{ v}}$ is surjective, then, for all
other ${\mathsf{ w}}\in{\mathbb{V}}$, \ ${{\sf image}\left[{\mathfrak{ g}}_{\mathsf{ w}}\right]}\subset Z({\mathcal{ B}})$.  If
$Z({\mathcal{ B}})=\{e\}$, then all other coefficients of ${\mathfrak{ g}}$ are trivial.

\item Suppose ${\mathsf{ v}}\neq{\mathsf{ w}}\in{\mathbb{V}}$ are such that
${\mathfrak{ g}}_{\mathsf{ v}} ={\mathfrak{ g}}_{\mathsf{ w}}$.  Then ${{\sf image}\left[{\mathfrak{ g}}_{\mathsf{ v}}\right]}$ is an abelian subgroup of
${\mathcal{ B}}$, and thus, $[{\mathcal{ B}},{\mathcal{ B}}] \subset \ker[{\mathfrak{ g}}_{\mathsf{ v}}]$.  Thus, if
$[{\mathcal{ B}},{\mathcal{ B}}]= {\mathcal{ B}}$, then ${\mathfrak{ g}}_{\mathsf{ v}}$ and ${\mathfrak{ g}}_{\mathsf{ w}}$ are trivial.

\item If  ${\mathcal{ B}}$ is simple but nonabelian, then only one
coefficient of ${\mathfrak{ G}}$ can be nontrivial; \ this coefficient is
an automorphism.  
\end{enumerate}
 \end{cor}
\bprf {\bf Part 1} and {\bf Part 2} are straightforward.  To see
{\bf Part 3}, note that $Z({\mathcal{ B}})$ is a
normal subgroup, so if ${\mathcal{ B}}$ is simple nonabelian, then $Z({\mathcal{ B}})=\{e\}$.  On
the other hand, any endomorphism of ${\mathcal{ B}}$ is either trivial or an
automorphism.  Hence, if ${\mathfrak{ G}}$ is nontrivial, it must have one
automorphic coefficient, and then, by {\bf Part 1}
 all other coefficients must be trivial.
{\tt \hrulefill $\Box$ }\end{list}  \medskip

\section{Structure Theory
\label{S:structure}}

  We now relate the structure of the group ${\mathcal{ B}}$ to the structure of
MCA on ${\mathcal{ B}}^{\mathbb{M}}$.  We review the structure theory of dynamical
systems in \S\ref{S:rel.dyn.sys} and group structure theory in
\S\ref{S:group.structure}.  In \S\ref{S:normal.decomposition}, we show
that, if ${\mathcal{ A}}$ is a fully characteristic subgroup of ${\mathcal{ B}}$, and
${\mathcal{ C}}={\mathcal{ B}}/{\mathcal{ A}}$, then the decomposition of ${\mathcal{ B}}$ into ${\mathcal{ A}}$
and ${\mathcal{ C}}$ yields a corresponding decomposition of MCA on
${\mathcal{ B}}^{\mathbb{M}}$.

\paragraph*{Notation:}  We will often decompose objects (eg. groups,
spaces, measures, functions) into factor and cofactor components.  We
will use three lexicographically consecutive letters to indicate,
respectively, the cofactor, product, and factor (eg. for groups:
${\mathcal{ A}}\hookrightarrow{\mathcal{ B}}\twoheadrightarrow{\mathcal{ C}}$; \ for measure spaces:
$({\mathbf{ Y}},{\mathcal{ Y}},\mu) = \left({\mathbf{ X}}\times{\mathbf{ Z}},
\ ,{\mathcal{ X}}\otimes{\mathcal{ Z}}, \ \lambda\otimes \nu\right)$; \
for dynamical systems, ${\sf G}={\sf F}\star{\sf H}$; \
for cellular automata, ${\mathfrak{ G}} = {\mathfrak{ F}}\star{\mathfrak{ H}}$, and for their local maps,
${\mathfrak{ g}} = {\mathfrak{ f}}\star{\mathfrak{ h}}$,  etc.).

\subsection{Relative and Nonhomogeneous CA
\label{S:relCA}
\label{S:rel.dyn.sys}}

  Let ${\mathbf{ X}}$ and ${\mathbf{ Z}}$ be a topological spaces.
 A topological  {\bf ${\mathbf{ Z}}$-relative dynamical
 system} \cite{Furstenberg,Petersen} on ${\mathbf{ X}}$ is a continuous  map
 ${\sf F}:{\mathbf{ X}} \times {\mathbf{ Z}}
{{\longrightarrow}} {\mathbf{ X}}$.  We write the second argument of ${\sf F}$ as a subscript: \ 
 for $\left( x, z\right)\in{\mathbf{ X}}\times{\mathbf{ Z}}$, \ \
 ${\sf F}( x, z)$ is written as
 ``${\sf F}_ z( x)$''.  Thus, ${\sf F}$ is treated as a
 ${\mathbf{ Z}}$-parameterized family of {\bf fibre maps}
 $\{{\sf F}_ z:{\mathbf{ X}}\!\longrightarrow{\mathbf{ X}}\!\}_{ z\in{\mathbf{ Z}}}$.  Let
 $\ensuremath{{\mathcal{ M}}\left[{\mathbf{ X}}\right] }$ be the set of Borel probability measures on
 ${\mathbf{ X}}$; \ if $\lambda\in\ensuremath{{\mathcal{ M}}\left[{\mathbf{ X}}\right] }$, then ${\sf F}$
 is {\bf $\lambda$-preserving} if ${\sf F}_ z (\lambda) \ = \
 \lambda$ for all $ z\in{\mathbf{ Z}}$.

  If ${\sf H}:{\mathbf{ Z}}\!\longrightarrow{\mathbf{ Z}}\!$ \ is a topological dynamical system,
 then the {\bf skew product} of ${\sf F}$ and ${\sf H}$ is the
 topological dynamical system ${\sf G} = {\sf F} \star {\sf H}$
 on ${\mathbf{ Y}} = {\mathbf{ X}} \times {\mathbf{ Z}}$ defined: \
$\displaystyle {\sf G}\left(  x, z\right) \ = \ \left( {\sf F}_ z( x),\  {\sf H}( z) \right)$.
  Now, suppose ${\mathbf{ X}} = {\mathcal{ A}}^{\mathbb{M}}$ and ${\mathbf{ Z}} = {\mathcal{ C}}^{\mathbb{M}}$, where
${\mathcal{ A}}$ and ${\mathcal{ C}}$ are finite sets. 
If ${\mathcal{ B}} = {\mathcal{ A}} \times {\mathcal{ C}}$, \  then there is a
 natural bijection ${\mathcal{ A}}^{\mathbb{M}} \times
{\mathcal{ C}}^{\mathbb{M}} \cong {\mathcal{ B}}^{\mathbb{M}}$.  A {\bf ${\mathcal{ C}}$-relative cellular automaton}
(RCA) on ${\mathcal{ A}}^{\mathbb{M}}$ is a continuous map ${\mathfrak{ F}}:{\mathcal{ A}}^{\mathbb{M}} \times
{\mathcal{ C}}^{\mathbb{M}} {{\longrightarrow}} {\mathcal{ A}}^{\mathbb{M}}$ which commutes with all shift maps:
$ {{{\boldsymbol{\sigma}}}^{{\mathsf{ m}}}} \circ {\mathfrak{ F}} = {\mathfrak{ F}} \circ
 {{{\boldsymbol{\sigma}}}^{{\mathsf{ m}}}} $, \ for all  ${\mathsf{ m}}\in{\mathbb{M}}$.  Like an ordinary CA,
${\mathfrak{ F}}$ is determined by a {\bf local map} ${\mathfrak{ f}}: {\mathcal{ A}}^{\mathbb{V}} \times
{\mathcal{ C}}^{\mathbb{V}} {{\longrightarrow}} {\mathcal{ A}}$, where ${\mathbb{V}}\subset{\mathbb{M}}$ is finite,
so that, for all $({\mathbf{ a}},{\mathbf{ c}}) \in {\mathcal{ B}}^{\mathbb{M}}$,  and ${\mathsf{ m}}\in{\mathbb{M}}$, \ \
${\mathfrak{ F}}({\mathbf{ a}},{\mathbf{ c}})_{\mathsf{ m}} \ = \ {\mathfrak{ f}}\left({\mathbf{ a}}\raisebox{-0.3em}{$\left|_{{\mathsf{ m}}+{\mathbb{V}}}\right.$},\
{\mathbf{ c}}\raisebox{-0.3em}{$\left|_{{\mathsf{ m}}+{\mathbb{V}}}\right.$}\right)$.  For any ${\mathbf{ c}}\in{\mathcal{ C}}^{\mathbb{V}}$, the {\bf local fibre map} ${\mathfrak{ f}}_{\mathbf{ c}}:{\mathcal{ A}}^{\mathbb{V}}{{\longrightarrow}}{\mathcal{ A}}$ is defined by
${\mathfrak{ f}}_{\mathbf{ c}}({\mathbf{ a}}) = {\mathfrak{ f}}({\mathbf{ c}},{\mathbf{ a}})$.  If ${\mathcal{ A}}$ is a
group and ${\mathfrak{ f}}_{\mathbf{ c}}$ is a product of affine endomorphisms for
every ${\mathbf{ c}} \in {\mathcal{ C}}^{\mathbb{V}}$, then ${\mathfrak{ F}}$ is called a {\bf
multiplicative} relative cellular automaton (MRCA).

  If ${\mathfrak{ H}}:{\mathcal{ C}}^{\mathbb{M}}\!\longrightarrow{\mathcal{ C}}^{\mathbb{M}}\!$ is a CA with local map
${\mathfrak{ h}}:{\mathcal{ C}}^{\mathbb{V}}{{\longrightarrow}}{\mathcal{ C}}$, then the {\bf skew product} ${\mathfrak{ F}} \star
{\mathfrak{ H}}$ is a CA on ${\mathcal{ B}}^{\mathbb{M}}$,
with local map
${\mathfrak{ g}}:{\mathcal{ B}}^{\mathbb{V}}\cong{\mathcal{ A}}^{\mathbb{V}}\times{\mathcal{ C}}^{\mathbb{V}}{{\longrightarrow}}{\mathcal{ B}}$ defined:
\ \ $\displaystyle  {\mathfrak{ g}}\left({\mathbf{ a}}, \ {\mathbf{ c}}\right) \ = \ \left( {\mathfrak{ f}}_{\mathbf{ c}}
 \left({\mathbf{ a}}\right), \ {\mathfrak{ h}}\left({\mathbf{ c}}\right)\right)$.
 
  A {\bf nonhomogeneous cellular automaton} (NHCA) is a continuous map
${\mathfrak{ G}}:{\mathcal{ B}}^{\mathbb{M}}\!\longrightarrow{\mathcal{ B}}^{\mathbb{M}}\!$ which does {\em not} necessarily commute with 
shift maps, but where there is some finite ${\mathbb{V}}\subset{\mathbb{M}}$, so that, 
for all ${\mathsf{ m}}\in{\mathbb{M}}$, there is a local map
${\mathfrak{ g}}_{{\mathsf{ m}}}:{\mathcal{ B}}^{({\mathsf{ m}}+{\mathbb{V}})}{{\longrightarrow}}{\mathcal{ B}}$ so that, 
$\forall {\mathbf{ b}}\in{\mathcal{ B}}^{\mathbb{M}}$, \ \ 
$\displaystyle {\mathfrak{ G}}({\mathbf{ b}})_{\mathsf{ m}} \ = \ {\mathfrak{ g}}_{{\mathsf{ m}}}\left({\mathbf{ b}}\raisebox{-0.3em}{$\left|_{({\mathsf{ m}}+{\mathbb{V}})}\right.$}\right)$. 
Thus, for example, any CA is an NHCA.  If
${\mathfrak{ F}}:{\mathcal{ A}}^{\mathbb{M}} \times {\mathcal{ C}}^{\mathbb{M}} {{\longrightarrow}} {\mathcal{ A}}^{\mathbb{M}}$ is an RCA,
then, for any ${\mathbf{ c}} \in {\mathcal{ C}}^{\mathbb{M}}$, the fibre map
${\mathfrak{ F}}_{\mathbf{ c}}:{\mathcal{ A}}^{\mathbb{M}}\!\longrightarrow{\mathcal{ A}}^{\mathbb{M}}\!$ is an NHCA.

\subsection{Group Structure Theory
\label{S:group.structure}}

 Let ${\mathcal{ B}}$ be a group. A subgroup ${\mathcal{ A}}\subset{\mathcal{ B}}$ is called {\bf fully
characteristic} \cite{WRScott} if, for every $\phi\in{{{\mathbf{ E}}{\mathbf{ n}}{\mathbf{ d}}_{} \left[{\mathcal{ B}}\right]}}$, we
have $\phi({\mathcal{ A}})\subset{\mathcal{ A}}$.  We indicate this:
``${\mathcal{ A}}\prec{\mathcal{ B}}$''.  For example, if $Z({\mathcal{ B}})$ is the center
of ${\mathcal{ B}}$, then  $Z({\mathcal{ B}})\prec {\mathcal{ B}}$. Observe that any fully characteristic
subgroup is normal.

 The main result of this section is:

\begin{thm}
\label{MCA.structure}   
Suppose that ${\mathcal{ A}}\prec{\mathcal{ B}}$, and ${\mathcal{ B}}/{\mathcal{ A}}={\mathcal{ C}}$.  If
${\mathfrak{ G}}:{\mathcal{ B}}^{\mathbb{M}}\!\longrightarrow{\mathcal{ B}}^{\mathbb{M}}\!$ is an MCA, then there is an
MCA ${\mathfrak{ H}}:{\mathcal{ C}}^{\mathbb{M}}\!\longrightarrow{\mathcal{ C}}^{\mathbb{M}}\!$ and an
MRCA ${\mathfrak{ F}}:{\mathcal{ A}}^{\mathbb{M}} \times {\mathcal{ C}}^{\mathbb{M}}{{\longrightarrow}}{\mathcal{ A}}^{\mathbb{M}}$ so that \
${\mathfrak{ G}} = {\mathfrak{ F}} \star {\mathfrak{ H}}$.
 \end{thm}

  We will prove this result in \S\ref{S:normal.decomposition}, and
also describes the structure of the local maps of ${\mathfrak{ H}}$ and ${\mathfrak{ F}}$
(see Proposition \ref{MCA.structure.prop}).  First we introduce the
relevant algebraic machinery.

\paragraph*{Semidirect Products and and Pseudoproducts:}

Suppose ${\mathcal{ A}}\subset{\mathcal{ B}}$ is a normal subgroup,
and ${\mathcal{ C}}={\mathcal{ B}}/{\mathcal{ A}}$, and
let $\pi:{\mathcal{ B}}\twoheadrightarrow{\mathcal{ C}}$ be the quotient map.
Let $\varsigma:{\mathcal{ C}} \rightarrowtail{\mathcal{ B}}$ be a section of $\pi$ ---that is,
for all $ c\in{\mathcal{ C}}$, \ $\pi\left(\varsigma( c)\right) =  c$. 
For any  $ a\in{\mathcal{ A}}$ and $ c\in{\mathcal{ C}}$, we define 
$ a \star  c :=  a \cdot \varsigma\left( c\right)$.
For every $b\in{\mathcal{ B}}$, there are unique  $ a\in{\mathcal{ A}}$ and
$ c\in{\mathcal{ C}}$ so that $b=  a \star  c$.
Thus, the map ${\mathcal{ A}} \times {\mathcal{ C}} \ni ( a, c)
\mapsto  a \star  c\in {\mathcal{ B}}$ is a bijection\footnote{...but
 generally not a homomorphism.}. 
 We  call ${\mathcal{ B}}$ 
a {\bf pseudoproduct} of ${\mathcal{ A}}$ and ${\mathcal{ C}}$, and
write: ``${\mathcal{ B}} = {\mathcal{ A}} \star {\mathcal{ C}}$''.

  If $ c\in{\mathcal{ C}}$, the {\bf conjugation automorphism}
$ c^{\ast}\,\!\in{{{\mathbf{ A}}{\mathbf{ u}}{\mathbf{ t}}_{} \left[{\mathcal{ A}}\right]}}$ is defined:
\[
  c^{\ast}\, a \ = \ \varsigma\left( c\right) \cdot  a \cdot\varsigma\left( c\right)^{-1}.
 \]
  Thus, multiplication using pseudoproduct notation satisfies the
equation:
\begin{equation}
 \left( a_1\star  c_1\right) \cdot
 \left( a_2\star  c_2\right) \ \ =  \ \ 
\left( a_1 \cdot \varsigma\left( c_1\right)\rule[-0.5em]{0em}{1em}\right) \cdot
 \left( a_2 \cdot \varsigma\left( c_2\right)\rule[-0.5em]{0em}{1em}\right) \ \ =  \ \ 
  a_1 \cdot \left( c_1^{\ast}\,  a_2\right) \cdot 
\left(\varsigma\left( c_1\right) \cdot \varsigma\left( c_2\right)\rule[-0.5em]{0em}{1em}\right).
\label{pseudodir.prod}
\end{equation}
  In general, $\varsigma\left( c_1\right) \cdot \varsigma\left( c_2\right)$ does not equal
$\varsigma\left( c_1\cdot  c_2\right)$; \ {\em this} is true only ${\mathcal{ B}}$ is a {\bf
semidirect product} of ${\mathcal{ A}}$ and ${\mathcal{ C}}$.  In this case, $\varsigma$
is an isomorphism from ${\mathcal{ C}}$ into an embedded
subgroup $\varsigma\left({\mathcal{ C}}\right)\subset{\mathcal{ B}}$, and (\ref{pseudodir.prod}) becomes:
\begin{equation}
 ( a_1\star  c_1) \cdot ( a_2\star  c_2) \ = \ 
 \left(\rule[-0.5em]{0em}{1em}  a_1 \cdot ( c_1^{\ast}\,  a_2)\right) \
 \star \ \left(\rule[-0.5em]{0em}{1em}  c_1 \cdot c_2\right).
\label{semidir.product}
\end{equation}
In this case, we write: ``${\mathcal{ B}} = {\mathcal{ A}} \rtimes {\mathcal{ C}}$''.  We can treat
${\mathcal{ C}}$ as embedded in ${\mathcal{ B}}$, so $\varsigma$ is
just the identity, and $ a\star c =  a\cdot c$.

  We call ${\mathcal{ B}}$  a {\bf polymorph} of ${\mathcal{ A}}$ if:
(1) ${\mathcal{ B}} = {\mathcal{ A}} \rtimes {\mathcal{ C}}$; \ 
(2) ${\mathcal{ A}}$ and $\varsigma\left({\mathcal{ C}}\right)$ are both fully characteristic
in ${\mathcal{ B}}$;  and (3) $ c^{\ast}\,\in  Z\left({{{\mathbf{ A}}{\mathbf{ u}}{\mathbf{ t}}_{} \left[{\mathcal{ A}}\right]}}\right)$,
for every $ c\in{\mathcal{ C}}$.

\medskip

  \medskip         \refstepcounter{thm} {\bf Example \thethm:}  \setcounter{enumi}{\thethm} \begin{list}{(\alph{enumii})}{\usecounter{enumii}} 			{\setlength{\leftmargin}{0em} 			\setlength{\rightmargin}{0em}}   
\item \label{X:Quat.pseudo}
 Let ${\mathcal{ B}}={\mathbf{ Q}}_8=\{\pm1,\ \pm {\mathbf{ i}},\ \pm {\mathbf{ j}},\ \pm {\mathbf{ k}}\} $ be the {\bf 
Quaternion Group}, defined by: ${\mathbf{ i}}^2 ={\mathbf{ j}}^2={\mathbf{ k}}^2= -1$,
and $q_1\cdot q_2 = q_3 = -q_2\cdot q_1$ for $(q_1,q_2,q_3)=({\mathbf{ i}},{\mathbf{ j}},{\mathbf{ k}})$ or any cyclic permutation thereof. 
Let ${\mathcal{ A}}=Z\left({\mathbf{ Q}}_8\right)=\{\pm1\}$; \  then
${\mathcal{ C}}={\mathbf{ Q}}_8/{\mathcal{ A}} \cong {{\mathbb{Z}}_{/2}}\oplus{{\mathbb{Z}}_{/2}}$.
Define homomorphism $\pi:{\mathbf{ Q}}_8 {{\longrightarrow}}{\mathcal{ C}}$
by $\pi(\pm 1)={\mathbf{ O}}:=(0,0)$, 
$\pi(\pm {\mathbf{ i}})={\mathbf{ I}}:=(1,0)$, 
$\pi(\pm {\mathbf{ j}})={\mathbf{ J}}:=(0,1)$, 
$\pi(\pm {\mathbf{ k}})={\mathbf{ K}}:=(1,1)$, with $\ker[\pi] = \{\pm 1\} = {\mathcal{ A}}$.  If $\varsigma:{\mathcal{ C}} {{\longrightarrow}} {\mathcal{ B}}$ is defined: $\varsigma\left({\mathbf{ O}}\right)=1$, 
$\varsigma\left({\mathbf{ I}}\right)={\mathbf{ i}}$, $\varsigma\left({\mathbf{ J}}\right)={\mathbf{ j}}$,  $\varsigma\left({\mathbf{ K}}\right)={\mathbf{ k}}$, then
$\varsigma$ induces a (non-semidirect) pseudoproduct structure ${\mathbf{ Q}}_8 = {\mathcal{ A}}\star {\mathcal{ C}}$.  In this case, ${\mathcal{ A}}=Z({\mathcal{ B}})$, and
multiplication satisfies the formula:
\[
  \left( a_1\star c_1\right) \cdot \left( a_2\star c_2\right)
\ \ = \ \ 
\left(\rule[-0.5em]{0em}{1em}  a_1\cdot a_2 \cdot \zeta\left( c_1, c_2\right)\right) \star \left( c_1\cdot c_2\right)
\hspace{2em}
\mbox{for all $ a_1, a_2\in{\mathcal{ A}}$, \ $ c_1, c_2\in{\mathcal{ C}}$.}
 \]
 Here, 
$\zeta\left( c_1, c_2\right) = 
\varsigma\left( c_1\cdot c_2\right)^{-1}\cdot
\varsigma( c_1)\cdot\varsigma\left( c_2\right) =
\sf sign\left[\varsigma\left( c_1\right)\cdot\varsigma\left( c_2\right)\right]$.
 For example, $\zeta({\mathbf{ I}},{\mathbf{ J}})=\sf sign[{\mathbf{ k}}]=
+1$, while $\zeta({\mathbf{ J}},{\mathbf{ I}})=\sf sign[-{\mathbf{ k}}] = -1$, and
 $\zeta\left({\mathbf{ O}}, c\right)=1$ for any $ c\in{\mathcal{ C}}$.

\item If $p$ is prime and ${\mathcal{ A}}=\left({{\mathbb{Z}}_{/p}},+\right)$ is the
(additive) cyclic group of order $p$, then ${{{\mathbf{ A}}{\mathbf{ u}}{\mathbf{ t}}_{} \left[{\mathcal{ A}}\right]}}$ is 
the (multiplicative) group $\left({{\mathbb{Z}}_{/p}}^\times,\cdot\right)$ of nonzero elements of the
field ${{\mathbb{Z}}_{/p}}$, acting on ${{\mathbb{Z}}_{/p}}$ by
multiplication, mod $p$.  The group $\left({{\mathbb{Z}}_{/p}}^\times, \cdot\right)$
is isomorphic to $\left({{\mathbb{Z}}_{/p-1}},+\right)$; \ thus, 
${\mathcal{ B}}= {{\mathbb{Z}}_{/p}}\rtimes {{\mathbb{Z}}_{/p}}^\times \cong {{\mathbb{Z}}_{/p}} \rtimes {{\mathbb{Z}}_{/p-1}}$ is a group of
order $p\cdot(p-1)$, and ${{\mathbb{Z}}_{/p}}$ is a characteristic subgroup. 
Since ${{\mathbb{Z}}_{/p}}^\times$ is itself abelian,
${\mathcal{ B}}$ is a polymorph of ${{\mathbb{Z}}_{/p}}$.

 If $q$ divides $p-1$, then there is a cyclic multiplicative subgroup
${\mathbf{ C}}_q\subset
{{\mathbb{Z}}_{/p}}^\times$ of order $q$. The semidirect product ${\mathbf{ D}}_{p;q}={{\mathbb{Z}}_{/p}}\rtimes {\mathbf{ C}}_q$,
is also a polymorph of ${{\mathbb{Z}}_{/p}}$. 
 For example, if $p=7$, then ${\mathbf{ C}}_3 = \{1,2,4\}$, and ${\mathbf{ D}}_{7;3}={{\mathbb{Z}}_{/7}}\rtimes{\mathbf{ C}}_3$ has cardinality 21.
  	\hrulefill\end{list}   			

\subsection{The Induced Decomposition
\label{S:normal.decomposition}}

  If ${\mathcal{ B}}={\mathcal{ A}}\star{\mathcal{ C}}$, then for any ${\mathbf{ a}}\in{\mathcal{ A}}^{\mathbb{M}}$ and
${\mathbf{ c}}\in{\mathcal{ C}}^{\mathbb{M}}$, define ${\mathbf{ b}}={\mathbf{ a}}\star{\mathbf{ c}}\in{\mathcal{ B}}^{\mathbb{M}}$ by
$b_{\mathsf{ m}} =  a_{\mathsf{ m}}\star c_{\mathsf{ m}}$ for all ${\mathsf{ m}}\in{\mathbb{M}}$;
we will thus identify ${\mathcal{ A}}^{\mathbb{M}}\times{\mathcal{ C}}^{\mathbb{M}}$ with ${\mathcal{ B}}^{\mathbb{M}}$.

  Suppose ${\mathcal{ B}} = {\mathcal{ A}}\rtimes{\mathcal{ C}}$, and let
${\mathfrak{ G}}:{\mathcal{ B}}^{\mathbb{Z}}\!\longrightarrow{\mathcal{ B}}^{\mathbb{Z}}\!$ be the {\bf nearest neighbour multiplication} CA,
Example (\ref{x0}).  If ${\mathbf{ b}}={\mathbf{ a}}\star{\mathbf{ c}}$, then
Moore \cite{MooreNLCA} noted that
\[
  {\mathfrak{ g}}(b_0,b_1) \ = \ 
\left( a_0\star  c_0\right) \cdot \left( a_1\star  c_1\right)
\ = \ 
\left( a_0\cdot  c_0^{\ast}\,  a_1\right) \star \left( c_0\cdot  c_1\right)
\ = \ 
{\mathfrak{ f}}_{\mathbf{ c}}\left( a_0, a_1\right) \star {\mathfrak{ h}}\left( c_0, c_1\right)
\]
where ${\mathfrak{ h}}:{\mathcal{ C}}^{\{0,1\}}{{\longrightarrow}}{\mathcal{ C}}$ is
defined ${\mathfrak{ h}}\left( c_0, c_1\right)=  c_0\cdot  c_1$, and 
${\mathfrak{ f}}:{\mathcal{ A}}^{\{0,1\}}\times{\mathcal{ C}}^{\{0,1\}}{{\longrightarrow}}{\mathcal{ A}}$ is
defined ${\mathfrak{ f}}_{\mathbf{ c}}\left( a_0, a_1\right)=  a_0\cdot  c_0^{\ast}\, a_1$.
Thus, ${\mathfrak{ G}} = {\mathfrak{ F}} \star {\mathfrak{ H}}$,
where ${\mathfrak{ H}}:{\mathcal{ C}}^{\mathbb{M}}\!\longrightarrow{\mathcal{ C}}^{\mathbb{M}}\!$ is the CA with local map ${\mathfrak{ h}}$,
and  ${\mathfrak{ F}}:{\mathcal{ A}}^{\mathbb{M}} \times {\mathcal{ C}}^{\mathbb{M}}{{\longrightarrow}}{\mathcal{ A}}^{\mathbb{M}}$ is the RCA
with local map ${\mathfrak{ f}}$.  In other words, the  decomposition
of ${\mathcal{ B}}={\mathcal{ A}}\rtimes{\mathcal{ C}}$ induces a  decomposition
of ${\mathfrak{ G}}$.  We now generalize this idea to arbitrary MCA.
 
\begin{lemma}\label{splitting.endomorphism.lemma}    Suppose ${\mathcal{ A}}\prec{\mathcal{ B}}$, \ ${\mathcal{ C}}={\mathcal{ B}}/{\mathcal{ A}}$, \ and
${\mathcal{ B}}={\mathcal{ A}}\star{\mathcal{ C}}$.  Let ${\mathfrak{ g}}\in{{{\mathbf{ E}}{\mathbf{ n}}{\mathbf{ d}}_{} \left[{\mathcal{ B}}\right]}}$.
\begin{enumerate}
\item  There exist
${\mathfrak{ f}}\in{{{\mathbf{ E}}{\mathbf{ n}}{\mathbf{ d}}_{} \left[{\mathcal{ A}}\right]}}$ and ${\mathfrak{ h}}\in{{{\mathbf{ E}}{\mathbf{ n}}{\mathbf{ d}}_{} \left[{\mathcal{ C}}\right]}}$
so that the following diagram commutes:
\[
\begin{array}{rcrcr}
{\mathcal{ A}} 		& \lefteqn{\hookrightarrow} &{\mathcal{ B}} 		&\lefteqn{\twoheadrightarrow} &{\mathcal{ C}} \\
{\mathfrak{ f}} \left\downarrow\rule[-0.5em]{0em}{1em}\right.&	&{\mathfrak{ g}} \left\downarrow\rule[-0.5em]{0em}{1em}\right.	&	  &{\mathfrak{ h}}\left\downarrow\rule[-0.5em]{0em}{1em}\right.\\
{\mathcal{ A}} 		& \lefteqn{\hookrightarrow} &{\mathcal{ B}} 		&\lefteqn{\twoheadrightarrow} &{\mathcal{ C}} 
\end{array}
\hspace{3em}\mbox{We indicate this: \ ``$\,{\mathfrak{ g}} = {\mathfrak{ f}} \star {\mathfrak{ h}}$.''}
\]
\item   Define ${\mathfrak{ g}}':{\mathcal{ C}}{{\longrightarrow}}{\mathcal{ A}}$ by ${\mathfrak{ g}}'( c)
\ = \ {\mathfrak{ g}}\left(\varsigma\left( c\right)\right)\cdot \varsigma\left({\mathfrak{ h}}\left( c\right)\right)^{-1}$.
If $ a\in{\mathcal{ A}}$ and $ c\in{\mathcal{ C}}$, \ then\ $\displaystyle
  {\mathfrak{ g}}(  a\star c) \ = \ 
\left(\rule[-0.5em]{0em}{1em}{\mathfrak{ f}}\left( a\right)\cdot {\mathfrak{ g}}'\left( c\right)\right) \ \star \
 {\mathfrak{ h}}\left( c\right)$.

\item If ${\mathcal{ B}}$ is a polymorph of ${\mathcal{ A}}$, then ${\mathfrak{ g}}'$ is trivial, 
so ${\mathfrak{ g}}(b) \ = \ {\mathfrak{ f}}\left( a\right)\star 
{\mathfrak{ h}}( c)$.

\item If ${\mathcal{ B}}$ is a polymorph of ${\mathcal{ A}}$ and ${\mathfrak{ f}} \in {{{\mathbf{ A}}{\mathbf{ u}}{\mathbf{ t}}_{} \left[{\mathcal{ A}}\right]}}$, then
${\mathfrak{ g}}\raisebox{-0.3em}{$\left|_{\varsigma\left({\mathcal{ C}}\right)}\right.$}={\mathbf{ Id}_{{}}}$ and  ${\mathfrak{ h}} = {{\mathbf{ Id}_{{_{{\mathcal{ C}}}}}}}$.

\end{enumerate}
 \end{lemma}
\bprf
  {\bf Part 1}: Define ${\mathfrak{ f}} = {\mathfrak{ g}}\raisebox{-0.3em}{$\left|_{{\mathcal{ A}}}\right.$}$.  Then
${\mathfrak{ f}}\in{{{\mathbf{ E}}{\mathbf{ n}}{\mathbf{ d}}_{} \left[{\mathcal{ A}}\right]}}$ because ${\mathcal{ A}}\prec{\mathcal{ B}}$.  Define
${\mathfrak{ h}}$ by: ${\mathfrak{ h}}(b\cdot {\mathcal{ A}}) = {\mathfrak{ g}}(b)\cdot{\mathcal{ A}}$ for any
coset $(b\cdot {\mathcal{ A}}) \in {\mathcal{ C}}$.  Then, for any
$(b_1\cdot{\mathcal{ A}})$ and $(b_2\cdot{\mathcal{ A}})$ in ${\mathcal{ C}}$, we have
${\mathfrak{ h}}(b_1{\mathcal{ A}}\cdot b_2{\mathcal{ A}}) \ = \ {\mathfrak{ h}}\left((b_1\cdot
b_2)\cdot {\mathcal{ A}}\right) \ = \ {\mathfrak{ g}}(b_1\cdot b_2)\cdot {\mathcal{ A}} \  = \
{\mathfrak{ g}}(b_1)\cdot{\mathfrak{ g}}(b_2)\cdot {\mathcal{ A}} \ = \
{\mathfrak{ g}}(b_1){\mathcal{ A}}\cdot{\mathfrak{ g}}(b_2){\mathcal{ A}} \ = \ {\mathfrak{ h}}(b_1{\mathcal{ A}})\cdot{\mathfrak{ h}}(b_2{\mathcal{ A}})$, so ${\mathfrak{ H}}$ is a endomorphism of ${\mathcal{ C}}$.
Clearly, ${\mathfrak{ h}}\circ\pi = \pi\circ{\mathfrak{ g}}$.

  {\bf Part 4}: 
If $ a_1, a_2\in{\mathcal{ A}}$
and $ d_1, d_2\in{\mathcal{ D}}=\varsigma\left({\mathcal{ C}}\right)$, then
by (\ref{semidir.product}),
\begin{eqnarray*}
\lefteqn{{\mathfrak{ f}}\left( a_1\right) \cdot\
  \left(\rule[-0.5em]{0em}{1em} {\mathfrak{ g}}\left( d_1\right)^{\ast}\,\circ{\mathfrak{ f}}\right)\left( a_2\right) \cdot 
\ \left(\rule[-0.5em]{0em}{1em} {\mathfrak{ g}}\left( d_1\right)\cdot {\mathfrak{ g}}\left( d_2\right)\right)}\\
& = &
\left(\rule[-0.5em]{0em}{1em} {\mathfrak{ g}}\left( a_1\right) \cdot {\mathfrak{ g}}\left( d_1\right)\right) \cdot 
\left(\rule[-0.5em]{0em}{1em} {\mathfrak{ g}}\left( a_2\right) \cdot {\mathfrak{ g}}\left( d_2\right)\right) 
\ \  = \ \  
{\mathfrak{ g}}\left(\rule[-0.5em]{0em}{1em}\left( a_1\cdot d_1\right) \cdot \left( a_2\cdot d_2\right)\right) \\
& = & 
{\mathfrak{ g}}\left(\rule[-0.5em]{0em}{1em}\left( a^{}_1\cdot d_1^{\ast}\,( a_2)\right)\cdot \left( d_1\cdot d_2\right)\right)
\ \ \  = \ \ \
 {\mathfrak{ f}}\left( a_1\right)\cdot\
   \left({\mathfrak{ f}}\circ d_1^{\ast}\,\right)( a_2) \ \cdot \ {\mathfrak{ g}}\left( d_1\right)\cdot
    {\mathfrak{ g}}\left( d_2\right).
\end{eqnarray*}
 Cancel  ${\mathfrak{ f}}\left( a_1\right)$ and ${\mathfrak{ g}}\left( d_1\right)
    {\mathfrak{ g}}\left( d_2\right)$, and note
that $ a_2$ is arbitrary to conclude:
\ ${\mathfrak{ g}}\left( d_1\right)^{\ast}\,\circ{\mathfrak{ f}} \ = \ {\mathfrak{ f}}\circ d_1^{\ast}\,$. \ 
Since $ d_1^{\ast}\,\in Z\left({{{\mathbf{ A}}{\mathbf{ u}}{\mathbf{ t}}_{} \left[{\mathcal{ A}}\right]}}\right)$, commute these
terms to get
\  ${\mathfrak{ g}}\left( d_1\right)^{\ast}\, \circ{\mathfrak{ f}} \ = \  d_1^{\ast}\,\circ{\mathfrak{ f}}$;
\ cancel ${\mathfrak{ f}}\in{{{\mathbf{ A}}{\mathbf{ u}}{\mathbf{ t}}_{} \left[{\mathcal{ A}}\right]}}$ to conclude that  
${\mathfrak{ g}}\left( d_1\right)^{\ast}\, \ = \  d_1^{\ast}\,$.
Now, ${\mathcal{ D}}$ is fully characteristic, so ${\mathfrak{ g}}\left( d_1\right)\in {\mathcal{ D}}$.
But the map ${\mathcal{ D}}\ni d\mapsto d^{\ast}\,\in{{{\mathbf{ A}}{\mathbf{ u}}{\mathbf{ t}}_{} \left[{\mathcal{ A}}\right]}}$ is
really just the inclusion map ${\mathcal{ D}}\cong{\mathcal{ C}}\hookrightarrow{{{\mathbf{ A}}{\mathbf{ u}}{\mathbf{ t}}_{} \left[{\mathcal{ A}}\right]}}$,
and therefore injective, so we conclude ${\mathfrak{ g}}\left( d_1\right) \ = \  d_1$.

  Thus, ${\mathfrak{ g}}\raisebox{-0.3em}{$\left|_{{\mathcal{ D}}}\right.$} = {{\mathbf{ Id}_{{_{{\mathcal{ D}}}}}}}$, so ${\mathfrak{ g}}\circ\varsigma=\varsigma$.
Since $\pi\circ\varsigma = {{\mathbf{ Id}_{{_{{\mathcal{ C}}}}}}}$, we conclude:
 ${\mathfrak{ h}} = {\mathfrak{ h}}\circ{{\mathbf{ Id}_{{_{{\mathcal{ C}}}}}}} =
 {\mathfrak{ h}}\circ\pi\circ\varsigma = \pi\circ{\mathfrak{ g}}\circ\varsigma
= \pi\circ\varsigma = {{\mathbf{ Id}_{{_{{\mathcal{ C}}}}}}}$.

{\bf Part 2} and {\bf Part 3} are straightforward.
{\tt \hrulefill $\Box$ }\end{list}  \medskip  
  \medskip         \refstepcounter{thm}                     \begin{list}{} 			{\setlength{\leftmargin}{1em} 			\setlength{\rightmargin}{0em}}                     \item {\bf Example \thethm:} \label{X:quat.map}
 Recall ${\mathbf{ Q}}_8={\mathcal{ A}}\star{\mathcal{ C}}$ from Example (\ref{X:Quat.pseudo}).
Define ${\mathfrak{ g}}_1,{\mathfrak{ g}}_2\in{{{\mathbf{ A}}{\mathbf{ u}}{\mathbf{ t}}_{} \left[{\mathbf{ Q}}_8\right]}}$ by Table \thethm.1.
Then ${\mathfrak{ f}}_1 = {\mathfrak{ f}}_2 = {{\mathbf{ Id}_{{_{{\mathcal{ A}}}}}}}$, while ${\mathfrak{ h}}_1$, ${\mathfrak{ h}}_2$,
${\mathfrak{ g}}_1'$, and ${\mathfrak{ g}}_2'$ are defined by Table \thethm.2.  	\hrulefill\end{list}   			  
\[
\begin{array}{ccc}
\begin{array}{|r||c|c|c|c|c|c|c|c|}
\hline
 & 1 & -1 & {\mathbf{ i}} & -{\mathbf{ i}} & {\mathbf{ j}} & -{\mathbf{ j}} & {\mathbf{ k}} & -{\mathbf{ k}} \\
\hline\hline
{\mathfrak{ g}}_1     & 1 & -1 & {\mathbf{ j}} & -{\mathbf{ j}} & {\mathbf{ k}} & -{\mathbf{ k}} & {\mathbf{ i}} & -{\mathbf{ i}} \\
\hline
{\mathfrak{ g}}_2     & 1 & -1 & -{\mathbf{ i}} & {\mathbf{ i}} & {\mathbf{ k}} & -{\mathbf{ k}} & {\mathbf{ j}} & -{\mathbf{ j}} \\
\hline
\end{array}
& \hspace{2em} &
\begin{array}{|r||c|c|c|c|}
\hline
    \rule[-0.5em]{0em}{1em}      & {\mathbf{ O}} & {\mathbf{ I}} & {\mathbf{ J}}  & {\mathbf{ K}}  \\
\hline\hline
{\mathfrak{ h}}_1  \rule[-0.5em]{0em}{1em}    & {\mathbf{ O}}  & {\mathbf{ J}} & {\mathbf{ K}} & {\mathbf{ I}} \\
\hline
{\mathfrak{ h}}_2  \rule[-0.5em]{0em}{1em}    & {\mathbf{ O}} &  {\mathbf{ I}} & {\mathbf{ K}} & {\mathbf{ J}}  \\
\hline\hline
{\mathfrak{ g}}_1'     & 1  & 1 & 1 & 1 \\
\hline
{\mathfrak{ g}}_2'    & 1  & -1 & 1 & 1 \\ 
\hline
\end{array}\\ \\
\mbox{Table \thethm.1} && \mbox{Table \thethm.2}
\end{array}
\]

\begin{prop}\label{MCA.structure.prop}  
The statement of Theorem \ref{MCA.structure} is true.
To be specific, if ${\mathfrak{ G}}$ has local map
\[
{\mathfrak{ g}}:{\mathcal{ B}}^{\mathbb{V}} \ni {\mathbf{ b}} \ \ \mapsto \ \
\left( g\cdot \prod_{i=1}^{I} {\mathfrak{ g}}_{i}\left(b_{ v[i]}\right)\right)
\in {\mathcal{ B}},
\hspace{1em}
\left(\mbox{where $b\in{\mathcal{ B}}$ and ${\mathfrak{ g}}_i\in{{{\mathbf{ E}}{\mathbf{ n}}{\mathbf{ d}}_{} \left[{\mathcal{ B}}\right]}}$, 
for all $i \in{\left[ 0..I \right]}$}\right)
\]
then the local maps ${\mathfrak{ h}}:{\mathcal{ C}}^{\mathbb{V}}{{\longrightarrow}}{\mathcal{ C}}$ and
and ${\mathfrak{ f}}:{\mathcal{ A}}^{\mathbb{V}}\times{\mathcal{ C}}^{\mathbb{V}}{{\longrightarrow}}{\mathcal{ A}}$ are defined as follows.
Fix a pseudoproduct representation ${\mathcal{ B}}={\mathcal{ A}}\star{\mathcal{ C}}$.  
Let $ g =  f \star  h$ for some $ f\in{\mathcal{ A}}$
and $ h\in{\mathcal{ C}}$.  For all $i\in{\left[ 0... I \right]}$,
let ${\mathfrak{ g}}_{i} = {\mathfrak{ f}}_{i} \star {\mathfrak{ h}}_{i}$, where
${\mathfrak{ f}}_{i}\in{{{\mathbf{ E}}{\mathbf{ n}}{\mathbf{ d}}_{} \left[{\mathcal{ A}}\right]}}$ and ${\mathfrak{ h}}_{i}\in{{{\mathbf{ E}}{\mathbf{ n}}{\mathbf{ d}}_{} \left[{\mathcal{ C}}\right]}}$,
as in Lemma \ref{splitting.endomorphism.lemma}.  Then:

\begin{enumerate}
\item \ $\displaystyle  {\mathfrak{ h}}({\mathbf{ c}}) \ \ =  \ \
 h\cdot \prod_{i=1}^{I} {\mathfrak{ h}}_{i}\left( c_{ v[i]}\right)$, \ 
and ${\mathfrak{ f}}$ is defined by expression (\ref{prpgf.2}) below. 

\hspace{-2.5em}In particular:
\item Suppose ${\mathcal{ B}}$ is a polymorph of ${\mathcal{ A}}$, and 
${\mathfrak{ f}}_i\in{{{\mathbf{ A}}{\mathbf{ u}}{\mathbf{ t}}_{} \left[{\mathcal{ B}}\right]}}$, $\forall i \in{\left[ 0..I \right]}$. 
Then 
$\displaystyle
{\mathfrak{ h}}\left( {\mathbf{ c}}\right) =   \prod_{i=0}^I 
  c_{ v[i]}$ \ and \
$\displaystyle
{\mathfrak{ f}}_{{\mathbf{ c}}}\left( {\mathbf{ a}}\right)  =   f \cdot  \prod_{i=0}^I 
{\mathfrak{ f}}^{i}_{{\mathbf{ c}}}\left(  a_{ v[i]}\right)$, \
where, \ $\forall i\geq0$, \ $\forall a\in{\mathcal{ A}}$, 
\ $\displaystyle {\mathfrak{ f}}^{i}_{{\mathbf{ c}}}( a)
 \ = \  h^{\ast}\, \,   c_{ v[0]}^{\ast}\, \,  c_{ v[1]}^{\ast}\,\, \ldots
\,   c_{ v[i-1]}^{\ast}\, \, {\mathfrak{ f}}_{i}( a)$.
\item Suppose ${\mathcal{ A}}\subset Z({\mathcal{ B}})$.  Treat ${\mathcal{ A}}$
as an additive group $({\mathcal{ A}},+)$.  Then 
${\mathfrak{ F}}_{{\mathbf{ c}}}({\mathbf{ a}}) = {\mathfrak{ L}}\left({\mathbf{ a}}\right)
 + {\mathfrak{ P}}\left({\mathbf{ c}}\right)$, where ${\mathfrak{ L}}:{\mathcal{ A}}^{\mathbb{M}}\!\longrightarrow{\mathcal{ A}}^{\mathbb{M}}\!$
is a linear cellular automaton with local map 
\begin{equation}
\label{xtlca.defn}
{\mathfrak{ l}}:{\mathcal{ A}}^{\mathbb{V}}\ni {\mathbf{ a}}
\ \ \mapsto \ \  \left(\sum_{i=0}^I {\mathfrak{ f}}_i\left( a_{ v[i]}\right) \right)
\in{\mathcal{ A}}
\end{equation} 
and ${\mathfrak{ P}}:{\mathcal{ C}}^{\mathbb{M}}{{\longrightarrow}}{\mathcal{ A}}^{\mathbb{M}}$ is a block map with local map 
${\mathfrak{ p}}:{\mathcal{ C}}^{\mathbb{V}}{{\longrightarrow}}{\mathcal{ A}}$ given by (\ref{gp.eqn}) below.
\end{enumerate}
 \end{prop}

  \medskip         \refstepcounter{thm} {\bf Example \thethm:}  \setcounter{enumi}{\thethm} \begin{list}{(\alph{enumii})}{\usecounter{enumii}} 			{\setlength{\leftmargin}{0em} 			\setlength{\rightmargin}{0em}}   
\item \label{relCAx0}
 Suppose ${\mathcal{ A}}\subset Z({\mathcal{ B}})$, as in {\bf Part 3} of Proposition
\ref{MCA.structure.prop}.  If
 $\displaystyle{\mathfrak{ g}}({\mathbf{ b}})= b_{{\mathsf{ v}}_1}^{n_1} b_{{\mathsf{ v}}_2}^{n_2}\cdots b_{{\mathsf{ v}}_J}^{n_J}$,\  then
 $\displaystyle {\mathfrak{ l}}\left({\mathbf{ a}}\right) \ = \ \sum_{{\mathsf{ v}}\in{\mathbb{V}}} \ell_{\mathsf{ v}}\cdot  a_{\mathsf{ v}}$,
where $\displaystyle\ell_{\mathsf{ v}}=\sum_{{\mathsf{ v}}_j={\mathsf{ v}}} n_j$ for each ${\mathsf{ v}}\in{\mathbb{V}}$.
Meanwhile, $\displaystyle{\mathfrak{ h}}({\mathbf{ c}})=  c_{{\mathsf{ v}}_1}^{n_1}
 c_{{\mathsf{ v}}_2}^{n_2}\cdots  c_{{\mathsf{ v}}_J}^{n_J}$.

\item \label{relCAx1}  Consider Example (\ref{x1}), with ${\mathcal{ B}} = {{\mathbb{Z}}_{/5}} \rtimes {{\mathbb{Z}}_{/4}}$.  In this case,
${\mathfrak{ h}}_i={{\mathbf{ Id}_{{_{{\mathcal{ C}}}}}}}$ and ${\mathfrak{ f}}_i={{\mathbf{ Id}_{{_{{\mathcal{ A}}}}}}}$ for all $i$.  Thus,
${\mathfrak{ H}}:\left({{\mathbb{Z}}_{/4}}\right)^{\mathbb{Z}}\!\longrightarrow\left({{\mathbb{Z}}_{/4}}\right)^{\mathbb{Z}}\!$ is the linear CA with
local map: \ 
${\mathfrak{ h}}\left( c_0, c_1, c_2\right) =  c_0 +  c_1 + c_2$.

 For any
 $ a\in{{\mathbb{Z}}_{/5}}$ and $ c\in{{\mathbb{Z}}_{/4}}$, \ \
$ c^{\ast}\, a = 2^{ c} \cdot  a$; \ 
thus,  {\bf Part 2} of Proposition \ref{MCA.structure.prop} implies that
${\mathfrak{ f}}:\left({{\mathbb{Z}}_{/5}}\right)^{{\left[ 0... 2 \right]}}\times\left({{\mathbb{Z}}_{/4}}\right)^{{\left[ 0... 2 \right]}}{{\longrightarrow}}{{\mathbb{Z}}_{/5}}$
is defined for all
$\left( a_0, a_1, a_2\right)\in({{\mathbb{Z}}_{/5}})^{{\left[ 0... 2 \right]}}$ and
$( c_0, c_1, c_2)\in({{\mathbb{Z}}_{/4}})^{{\left[ 0... 2 \right]}}$, by: \ 
$\displaystyle
 {\mathfrak{ f}}_{( c_0, c_1, c_2)} \left( a_0, a_1, a_2\right)
 \ = \  a_0 + 2^{ c_0}  a_1 + 2^{ c_0+ c_1}  a_2$.

\item \label{relCAx2} Consider Example (\ref{x2}), with ${\mathcal{ B}} = {{\mathbb{Z}}_{/5}} \rtimes {{\mathbb{Z}}_{/4}}$.
Now  ${\mathfrak{ H}}:\left({{\mathbb{Z}}_{/4}}\right)^{\mathbb{Z}}\!\longrightarrow\left({{\mathbb{Z}}_{/4}}\right)^{\mathbb{Z}}\!$ has local map
$\displaystyle
 {\mathfrak{ h}}\left( c_0, c_1, c_2\right) \ = \  c_0 + 3 c_1 +
4 c_2 \ \equiv \  c_0- c_1 \pmod{4}$.
Meanwhile, {\bf Part 2} of Proposition \ref{MCA.structure.prop} implies that
\begin{eqnarray*}
\lefteqn{  {\mathfrak{ f}}_{( c_0, c_1, c_2)} \left( a_0, a_1, a_2\right)}\\
 & = &
 a_2 + 2^{ c_2}  a_2  
+ 2^{2 c_2}  a_2 
+ 2^{3 c_2}  a_2  
+ 2^{4 c_2}  a_1 
+ 2^{4 c_2+ c_1}  a_1  
+ 2^{4 c_2+2 c_1}  a_1  
+ 2^{4 c_2+3 c_1}  a_0 \\
& = &
 \left( 1 +  2^{ c_2} +  2^{2 c_2} + 2^{3 c_2}\right)  a_2
\ + \ \left(2^{4 c_2} +  2^{4 c_2+ c_1} + 2^{4 c_2+2 c_1}\right)  a_1
\ + \ 2^{4 c_2+3 c_1}  a_0.
\end{eqnarray*}

\item
 Recall  ${\mathfrak{ g}}_1,{\mathfrak{ g}}_2\in{{{\mathbf{ A}}{\mathbf{ u}}{\mathbf{ t}}_{} \left[{\mathbf{ Q}}_8\right]}}$ from Example \ref{X:quat.map},
and define ${\mathfrak{ g}}:{\mathbf{ Q}}_8^{{\left[ 0... 2 \right]}}{{\longrightarrow}}{\mathbf{ Q}}_8$ by ${\mathfrak{ g}}(q_0,q_1,q_2)=
q_0\cdot {\mathfrak{ g}}_1(q_1) \cdot {\mathfrak{ g}}_2(q_2)$. 
Identify $\{\pm1\}={\mathcal{ A}}$ with $\left({{\mathbb{Z}}_{/2}},+\right)$, and identify
${\mathcal{ C}}$ with ${{\mathbb{Z}}_{/2}}\oplus{{\mathbb{Z}}_{/2}}$ as described
in Example (\ref{X:Quat.pseudo}).
Then by {\bf Part 3} of Proposition \ref{MCA.structure.prop}, 
\[
\begin{array}{rcccl}
{\mathfrak{ l}}\left( a_0, a_1, a_2\right)
   & = &
  a_0 + {\mathfrak{ f}}_1\left( a_1\right)+
      {\mathfrak{ f}}_2\left( a_2\right) & = &  a_0 +  a_1 +  a_2;
 \\ \\
{\mathfrak{ h}}\left( c_0, c_1, c_2\right)
   & = &
  c_0 + {\mathfrak{ h}}_1\left( c_1\right) +  {\mathfrak{ h}}_2\left( c_2\right) 
&=&
\displaystyle \left[ c_{0,1}\atop c_{0,2}\right]
+{\left[\begin{array}{ccccccccccccccccccccccccr} 0&1\\1&1 \end{array}\right]} \left[ c_{1,1}\atop c_{1,2}\right]
+{\left[\begin{array}{ccccccccccccccccccccccccr} 1&1\\0&1 \end{array}\right]} \left[ c_{2,1}\atop c_{2,2}\right].
\end{array}
\]
 (where $ c_i \ = \  \left[ c_{i,1}\atop c_{i,2}\right]
\ \in {{\mathbb{Z}}_{/2}}\oplus{{\mathbb{Z}}_{/2}}$, for $i=0,1,2$).
Also, applying (\ref{gp.eqn}) below,
\[
{\mathfrak{ p}}\left( c_0, c_1, c_2\right)
   \ = \ 
  {\mathfrak{ e}}\left( c_0, c_1, c_2\right)+  
 {\mathbf{ Id}_{{}}}'\left( c_0\right)+{\mathfrak{ g}}_1'\left( c_1\right)  + {\mathfrak{ g}}_2'\left( c_2\right)
   \ = \ 
  {\mathfrak{ e}}\left( c_0, c_1, c_2\right) +  {\mathfrak{ g}}_2'\left( c_2\right),
\]
where ${\mathfrak{ e}}\left( c_0, c_1, c_2\right) = 
{\left\{ \begin{array}{rcl}                                 0 &&\mbox{if} \ \
 \varsigma\left( c_0\right)\cdot
\varsigma\left( c_1\right)\cdot
\varsigma\left({\mathfrak{ h}}_2\left( c_2\right)\right)
\ = \ \varsigma\left[  c_0 \cdot  c_1 \cdot {\mathfrak{ h}}_2\left( c_2\right)\right]\\
1 && \mbox{if} \ \ \varsigma\left( c_0\right)\cdot
\varsigma\left( c_1\right)\cdot
\varsigma\left({\mathfrak{ h}}_2\left( c_2\right)\right)
\ = \ -\varsigma\left[  c_0 \cdot  c_1 \cdot {\mathfrak{ h}}_2\left( c_2\right)\right]
                                  \end{array}  \right.  }$.

\item  \label{X:quat.decomp} ${\mathfrak{ g}}:{\mathbf{ Q}}_8^{{\left[ 0... 3 \right]}}{{\longrightarrow}}{\mathbf{ Q}}_8$ by 
${\mathfrak{ g}}(q_0,q_1,q_2,q_3)= q_3\cdot q_0^3\cdot q_2^5 \cdot q_1^{-1}$.
Then, in additive notation, 
${\mathfrak{ l}}\left( a_0, a_1, a_2, a_3\right)
\ = \  a_0 +  a_1 +   a_2 +  a_3$
and ${\mathfrak{ h}}\left( c_0, c_1, c_2, c_3\right)
\ = \   c_0 +  c_1 +   c_2 +  c_3 \pmod{2}$.\hrulefill
 	\hrulefill\end{list}   			

\begin{list}{} 			{\setlength{\leftmargin}{1em} 			\setlength{\rightmargin}{0em}}                         \item {\bf \hspace{-1em}  Proof of Proposition \ref{MCA.structure.prop}: \ \ } 
For $i\in{\left[ 0... I \right]}$, define
${\mathfrak{ f}}^{i}:{\mathcal{ C}}^{\mathbb{V}}\times{\mathcal{ A}}{{\longrightarrow}}{\mathcal{ A}}$ by
\begin{equation}
 \label{prpgf.1}
{\mathfrak{ f}}^{i}_{{\mathbf{ c}}}( a)
\ \ =  \ \  
\displaystyle \left( h \cdot \prod_{j=0}^{i-1} {\mathfrak{ h}}_j\left( c_{ v[j]}\right)\right)^{\ast}\,
\left(\rule[-0.5em]{0em}{1em} {\mathfrak{ f}}_{i}  \left( a\right)\cdot {\mathfrak{ g}}_i'\left( c_{ v[i]}\right)\right);
\end{equation}

for example, $\displaystyle {\mathfrak{ f}}^{0}_{{\mathbf{ c}}}( a)
 \ = \
   h^{\ast}\, \left(\rule[-0.5em]{0em}{1em} {\mathfrak{ f}}_{0}( a) \cdot  {\mathfrak{ g}}'_0\left( c_{ v[0]}\right)\right)$.  
Next, define ${\mathfrak{ e}}:{\mathcal{ C}}^{\mathbb{V}} {{\longrightarrow}} {\mathcal{ A}}$ by:
\begin{equation}
\label{tln.defn}
 {\mathfrak{ e}}({\mathbf{ c}}) \ = \ \left(\varsigma\left( h\right) \cdot 
 \prod_{i=0}^I \varsigma\left({\mathfrak{ h}}_{i}\left( c_{ v[i]}\right)\right) \right)
\cdot 
 \varsigma\left( h\cdot  \prod_{i=0}^I {\mathfrak{ h}}_{i}\left( c_{ v[i]}\right) \right)^{-1},
\end{equation}
  and define ${\mathfrak{ f}}:{\mathcal{ A}}^{\mathbb{V}} \times {\mathcal{ C}}^{\mathbb{V}}{{\longrightarrow}}{\mathcal{ A}}$ by: \
\begin{equation}
\label{prpgf.2}
\displaystyle {\mathfrak{ f}}_{{\mathbf{ c}}}\left( {\mathbf{ a}}\right) \ = \  f \cdot \left( \prod_{i=0}^I 
{\mathfrak{ f}}^{i}_{{\mathbf{ c}}}\left(  a_{ v[i]}\right) \right) \cdot {\mathfrak{ e}}({\mathbf{ c}}).
\end{equation} 

To prove {\bf Part 1}, we must show, for any
 ${\mathbf{ a}}\in{\mathcal{ A}}^{\mathbb{V}}$ and ${\mathbf{ c}}\in{\mathcal{ C}}^{\mathbb{V}}$, that $\displaystyle
 {\mathfrak{ g}}({\mathbf{ a}}\star{\mathbf{ c}}) \ = \ {\mathfrak{ f}}_{{\mathbf{ c}}}({\mathbf{ a}})
 \star {\mathfrak{ h}}({\mathbf{ c}})$. 

To see this, let
 ${\mathbf{ c}} = {\left[ c_{\mathsf{ v}}  |_{{\mathsf{ v}}\in{\mathbb{V}}}^{} \right]}
\in {\mathcal{ C}}^{\mathbb{V}}$ and ${\mathbf{ a}} = {\left[ a_{\mathsf{ v}}  |_{{\mathsf{ v}}\in{\mathbb{V}}}^{} \right]} \in {\mathcal{ A}}^{\mathbb{V}}$.   Then
\[
 {\mathfrak{ g}}\left({\mathbf{ a}}\star{\mathbf{ c}}\right) \ \ \ = \ \ \ 
 g\cdot
\prod_{i=0}^I {\mathfrak{ g}}_{i}\left(\rule[-0.5em]{0em}{1em} a_{ v[i]}\cdot \varsigma\left( c_{ v[i]}\right)\right) 
\  \ \ = \ \ \ 
 f\cdot \varsigma\left( h\right) \cdot \prod_{i=0}^I  {\mathfrak{ g}}_{i}\left( a_{ v[i]}\right)\cdot 
{\mathfrak{ g}}_{i}\left(\varsigma\left( c_{ v[i]}\right)\right).
\]
In the case $I=1$, this becomes:
\begin{eqnarray*}
\lefteqn{ {\mathfrak{ g}}\left({\mathbf{ a}}\star{\mathbf{ c}}\right) \ \ = \ \ 
  f\cdot \varsigma\left( h\right) \ \cdot\  
 {\mathfrak{ g}}_{0}\left( a_{ v[0]}\right)\cdot 
{\mathfrak{ g}}_{0}\left(\varsigma\left( c_{ v[0]}\right)\right) \ \cdot \ 
 {\mathfrak{ g}}_{1}\left( a_{ v[1]}\right)\cdot 
{\mathfrak{ g}}_{1}\left(\varsigma\left( c_{ v[1]}\right)\right)}
\\
& = &
 f\cdot \varsigma\left( h\right) \ \cdot\  
 {\mathfrak{ f}}_{0}\left( a_{ v[0]}\right)\cdot {\mathfrak{ g}}_0'\left( c_{ v[0]}\right) \cdot
\varsigma\left({\mathfrak{ h}}_{0}\left( c_{ v[0]}\right)\right) \ \cdot \ 
 {\mathfrak{ f}}_{1}\left( a_{ v[1]}\right)\cdot  {\mathfrak{ g}}_1'\left( c_{ v[1]}\right) \cdot
\varsigma\left({\mathfrak{ h}}_{1}\left( c_{ v[1]}\right)\right)\\
& = &
  f\cdot \varsigma\left( h\right) \ \cdot\  
 {\mathfrak{ f}}_{0}\left( a_{ v[0]}\right)\cdot {\mathfrak{ g}}_0'\left( c_{ v[0]}\right) 
\ \ \cdot \ \ {\mathfrak{ h}}_{0}( c_{ v[0]})^{\ast}\,\!
\left(\rule[-0.5em]{0em}{1em} {\mathfrak{ f}}_{1}\left( a_{ v[1]}\right)  {\mathfrak{ g}}_1'\left( c_{ v[1]}\right) \right)
\\&&\hspace{20em}
\ \  \ \cdot  \
\varsigma\left({\mathfrak{ h}}_{0}\left( c_{ v[0]}\right)\right) \cdot  
\varsigma\left({\mathfrak{ h}}_{1}\left( c_{ v[1]}\right)\right)\\
& = &
  f\ \cdot \   h^{\ast}\,  
 \left(\rule[-0.5em]{0em}{1em}{\mathfrak{ f}}_{0}\left( a_{ v[0]}\right)\cdot {\mathfrak{ g}}_0'\left( c_{ v[0]}\right)\right) 
\ \ \cdot \ \
 \left(  h \cdot {\mathfrak{ h}}_{0}\left( c_{ v[0]}\right)\right)^{\ast}\,
 \left(\rule[-0.5em]{0em}{1em} {\mathfrak{ f}}_{1}\left( a_{ v[1]}\right) \cdot {\mathfrak{ g}}_1'\left( c_{ v[1]}\right) \right)
\\&&\hspace{20em}
  \ \cdot \ 
 \varsigma\left( h\right) \cdot
\varsigma\left({\mathfrak{ h}}_{0}\left( c_{ v[0]}\right)\right) \cdot  
\varsigma\left({\mathfrak{ h}}_{1}\left( c_{ v[1]}\right)\right)\\
& = &
  f\cdot {\mathfrak{ f}}^{0}_{{\mathbf{ c}}}\left( a_{ v[0]}\right) \cdot
  {\mathfrak{ f}}^{1}_{{\mathbf{ c}}}\left(  a_{ v[1]}\right) 
 \ \cdot \ 
{\mathfrak{ e}}({\mathbf{ c}}) \cdot
 \varsigma\left(\rule[-0.5em]{0em}{1em}  h \cdot
{\mathfrak{ h}}_{0}\left( c_{ v[0]}\right) \cdot  
{\mathfrak{ h}}_{1}\left( c_{ v[1]}\right) \right) \\
& = &
{\mathfrak{ f}}_{{\mathbf{ c}}}({\mathbf{ a}})\cdot \varsigma\left( {\mathfrak{ h}}({\mathbf{ c}}) \right) 
\ \ \  =  \ \ \ 
{\mathfrak{ f}}_{{\mathbf{ c}}}({\mathbf{ a}}) \star {\mathfrak{ h}}({\mathbf{ c}}). 
\end{eqnarray*}
 A similar argument clearly works for $I\geq2$.

  It remains to show that ${\mathfrak{ e}}({\mathbf{ c}}) \in {\mathcal{ A}}$, which is
equivalent to showing that $\pi\left({\mathfrak{ e}}({\mathbf{ c}})\right) =  e_{{\mathcal{ C}}}$, where
$\pi:{\mathcal{ B}}\twoheadrightarrow{\mathcal{ C}}$ is the quotient map and $ e_{{\mathcal{ C}}}\in{\mathcal{ C}}$ is the
identity.  But $\\ \displaystyle \pi\left(\rule[-0.5em]{0em}{1em}{\mathfrak{ e}}({\mathbf{ c}})\right) \ = \
\pi\left[\left( \varsigma\left( h\right) \cdot 
 \prod_{i=0}^I \varsigma\left({\mathfrak{ h}}_{i}\left( c_{ v[i]}\right)\right) \right)
\cdot 
\varsigma\left(
  h\cdot  \prod_{i=0}^I {\mathfrak{ h}}_{i}\left( c_{ v[i]}\right)\right)^{-1}
 \right] $

$\displaystyle \hspace{4.5em} = \
\pi\left(\rule[-0.5em]{0em}{1em}\varsigma\left( h\right)\right) \cdot 
 \prod_{i=0}^I \pi\left(\rule[-0.5em]{0em}{1em}\varsigma\left({\mathfrak{ h}}_{i}\left( c_{ v[i]}\right)\right)\right)
\cdot
\pi\left(
\varsigma\left(  h\cdot  \prod_{i=0}^I {\mathfrak{ h}}_{i}\left( c_{ v[i]}\right) \right)\right)^{-1}  $

$\displaystyle \hspace{4.5em} = \
 h \cdot  \prod_{i=0}^I {\mathfrak{ h}}_{i}\left( c_{ v[i]}\right)
\cdot
\left(  h\cdot  \prod_{i=0}^I {\mathfrak{ h}}_{i}\left( c_{ v[i]}\right) \right)^{-1}
\ \ \  = \ \  \  e_{{\mathcal{ C}}}$.

{{\bf Part 2}:} ${\mathfrak{ g}}'_i$ and ${\mathfrak{ e}}$ are trivial, and
{\bf Part 4} of Lemma \ref{splitting.endomorphism.lemma}
 implies that ${\mathfrak{ h}}_i={{\mathbf{ Id}_{{_{{\mathcal{ C}}}}}}}$, \ $\forall i\in{\left[ 0... I \right]}$. 

{{\bf Part 3}:} All the conjugation automorphisms are trivial,
so expression (\ref{prpgf.1}) (written additively) simplifies to
\ ${\mathfrak{ f}}^{i}_{{\mathbf{ c}}}( a)
 \ = \ {\mathfrak{ f}}_{i}\left( a\right)
  \ + \  {\mathfrak{ g}}_i'\left({\mathbf{ c}}\right)$, for all $i\in{\left[ 0... I \right]}$.  Substitute
this into (\ref{prpgf.2}), to get \
$\displaystyle {\mathfrak{ f}}_{{\mathbf{ c}}}({\mathbf{ a}}) \   = \ \ \ 
  f \ + \ \sum_{i=0}^I \left(\rule[-0.5em]{0em}{1em} {\mathfrak{ f}}_{i}\left( a\right)
  \ + \  {\mathfrak{ g}}_i'\left({\mathbf{ c}}\right) \right) \ + \ {\mathfrak{ e}}({\mathbf{ c}})
\  = \ \ \ 
{\mathfrak{ l}}\left({\mathbf{ a}}\right)
 \ + \  {\mathfrak{ p}}\left({\mathbf{ c}}\right)$, \ \ 
 where \
\begin{equation}
\label{gp.eqn}
 {\mathfrak{ p}}\left({\mathbf{ c}}\right) \ \ = \  \  f \ + \  {\mathfrak{ e}}\left({\mathbf{ c}}\right) \ + \ 
\sum_{i=0}^I {\mathfrak{ g}}'_i\left({\mathbf{ b}}_{ v[i]}\right), 
\end{equation} 
and ${\mathfrak{ l}}\left({\mathbf{ a}}\right)$ is as in (\ref{xtlca.defn}).
{\tt \hrulefill $\Box$ }\end{list}  \medskip

\paragraph*{Application to Nilpotent groups:}{\label{upper.central.MCA}

A {\bf fully characteristic series} is an ascending chain of subgroups:
\begin{equation}
\label{composition.series}
\label{upper.central.series}
 \{e\} \ = \ {\mathcal{ Z}}_0 \prec {\mathcal{ Z}}_1 \prec \ldots 
\prec {\mathcal{ Z}}_K = {\mathcal{ B}}. 
\end{equation}
  where each is fully characteristic in the next.
 For example, the {\bf upper central series} of ${\mathcal{ B}}$ is the series
(\ref{composition.series}), where ${\mathcal{ Z}}_1 = Z({\mathcal{ B}})$, and for each
$k\geq1$, \ ${\mathcal{ Z}}_{k+1}$ is the complete preimage in ${\mathcal{ B}}$ of $Z\left(
{\mathcal{ C}}_{k}\right)$ under the quotient map ${\mathcal{ B}}\twoheadrightarrow
{\mathcal{ C}}_{k}:={\mathcal{ B}}/{\mathcal{ Z}}_{k}$, until we reach $K>0$ so that ${\mathcal{ Z}}_K =
{\mathcal{ Z}}_{K+1} = {\mathcal{ Z}}_{K+2} = \ldots$.  Thus, for all $k\in{\left[ 1... K \right]}$, the
factor groups ${\mathcal{ Q}}_k = {\mathcal{ Z}}_{k}/{\mathcal{ Z}}_{k-1}
\ \cong \ Z\left({\mathcal{ B}}/{\mathcal{ Z}}_{k-1}\right)$
are abelian (but ${\mathcal{ C}}_K = {\mathcal{ B}}/{\mathcal{ Z}}_{K}$ is not).
In general, ${\mathcal{ Z}}_K \neq {\mathcal{ B}}$; \ if they
are equal, then ${\mathcal{ B}}$ is called {\bf nilpotent}, and ${\mathcal{ C}}$ is trivial.

  \medskip         \refstepcounter{thm}                     \begin{list}{} 			{\setlength{\leftmargin}{1em} 			\setlength{\rightmargin}{0em}}                     \item {\bf Example \thethm:} \label{X:quaternions}
 Let ${\mathcal{ B}}={\mathbf{ Q}}_8$ from Example \ref{X:Quat.pseudo}.
 Then ${\mathcal{ Z}}_1=Z\left({\mathbf{ Q}}_8\right)=\{\pm1\}$, and
${\mathbf{ Q}}_8/{\mathcal{ Z}}_1 \cong {{\mathbb{Z}}_{/2}}\oplus{{\mathbb{Z}}_{/2}}$ is abelian,
so that ${\mathcal{ Z}}_2 = {\mathbf{ Q}}_8$.  Thus, ${\mathbf{ Q}}_8$ is nilpotent,
with upper central series:  \ $\displaystyle \{1\} \prec \{\pm1\} \prec {\mathbf{ Q}}_8$.
  	\hrulefill\end{list}   			  

 We can apply Theorem \ref{MCA.structure} recursively to a totally
 characteristic series like (\ref{composition.series}).  Let
 ${\mathcal{ A}}_1={\mathcal{ Z}}_1$, and ${\mathcal{ C}}_1 = {\mathcal{ B}}/{\mathcal{ A}}_1$, and write
 ${\mathcal{ B}}={\mathcal{ A}}_1\star{\mathcal{ C}}_1$, so that ${\mathfrak{ G}} = {\mathfrak{ F}}_{1}
 \star {\mathfrak{ H}}_1$, where ${\mathfrak{ H}}_1:{\mathcal{ C}}_1^{\mathbb{M}}\!\longrightarrow{\mathcal{ C}}_1^{\mathbb{M}}\!$ and
 ${\mathfrak{ F}}_{1}:{\mathcal{ A}}_1^{\mathbb{M}}\times{\mathcal{ C}}_1^{\mathbb{M}} {{\longrightarrow}} {\mathcal{ A}}_1^{\mathbb{M}}$ are
 multiplicative.  The series (\ref{composition.series}) induces a
 fully characteristic series
\begin{equation}
\label{factored.composition.series}
 \{e\} \ = \ {\mathcal{ Z}}_{_{1/1}} \prec {\mathcal{ Z}}_{_{2/1}} \prec \ldots 
\prec {\mathcal{ Z}}_{_{K/1}} = {\mathcal{ C}}_1,
\end{equation}
  where for each $k\in{\left[ 1... K \right]}$, we let ${\mathcal{ Z}}_{_{k/1}} =
{\mathcal{ Z}}_k/{\mathcal{ A}}_1\subset{\mathcal{ C}}_1$.  Now let ${\mathcal{ A}}_2={\mathcal{ Z}}_{_{2/1}}$ and
 ${\mathcal{ C}}_2 = {\mathcal{ C}}_1/{\mathcal{ A}}_2
\cong {\mathcal{ B}}/{\mathcal{ Z}}_2$, and write ${\mathcal{ C}}_1 = {\mathcal{ A}}_2\star{\mathcal{ C}}_2$,
so that ${\mathfrak{ H}}_1={\mathfrak{ F}}_2\star {\mathfrak{ H}}_{2}$, where
 ${\mathfrak{ F}}_{2}:{\mathcal{ A}}_{2}^{\mathbb{M}}\times
{\mathcal{ C}}_2^{\mathbb{M}}{{\longrightarrow}} {\mathcal{ A}}_{2}^{\mathbb{M}}$ and
${\mathfrak{ H}}_2:{\mathcal{ C}}_2^{\mathbb{M}}\!\longrightarrow{\mathcal{ C}}_2^{\mathbb{M}}\!$.  Proceed inductively.  In
particular, if ${\mathcal{ B}}$ is nilpotent, apply this to the upper central
series to obtain a decomposition:\ ${\mathfrak{ G}} =
{\mathfrak{ F}}_1\star\left({\mathfrak{ F}}_2\star\left[\ldots\star\left({\mathfrak{ F}}_{K-1}
\star{\mathfrak{ H}}\right)\ldots\right]\right)$, where, for all $k\in{\left[ 1... K \right)}$, \ 
 ${\mathfrak{ F}}_k:{\mathcal{ A}}_k^{\mathbb{M}}\times{\mathcal{ C}}_{k}^{\mathbb{M}}{{\longrightarrow}}{\mathcal{ A}}_k^{\mathbb{M}}$
is an affine RCA, while ${\mathfrak{ H}}:{\mathcal{ A}}_K^{\mathbb{M}}\!\longrightarrow{\mathcal{ A}}_K^{\mathbb{M}}\!$ is an affine
CA on the final abelian factor ${\mathcal{ A}}_K = {\mathcal{ B}}/{\mathcal{ Z}}_{K-1}$.
\section{Entropy \label{S:entropy}}

  Throughout this section, ${\mathbb{M}}={\mathbb{Z}}$ or ${\mathbb{N}}$, and
${\mathbb{V}}={\left[  V_0... V_1 \right]}\subset{\mathbb{M}}$.
Let $ L = -\min\{ V_0,0\}$ and $ R = \max\{0, V_1\}$,
and let $ V= R +  L$.   Let ${\mathcal{ B}}$ be a finite set, with $ B={{\sf card}\left[{\mathcal{ B}}\right]}$.
If ${\mathbf{ c}}\in {\mathcal{ B}}^{{\left[  J... K \right)}}$ and ${\left\langle {\mathbf{ c}} \right\rangle } =
{\left\{ {\mathbf{ b}}\in{\mathcal{ B}}^{\mathbb{Z}} \; ; \; {\mathbf{ b}}\raisebox{-0.3em}{$\left|_{{\left[  J... K \right)}}\right.$} = {\mathbf{ c}} \right\} }$ is the
corresponding {\bf cylinder set}, we say ${\left\langle {\mathbf{ c}} \right\rangle }$ is a cylinder set
of {\bf length $\ell= K- J$}.  Let $\eta_{_{{\mathcal{ B}}}}{}$ be the
uniformly distributed Bernoulli measure on ${\mathcal{ B}}^{\mathbb{Z}}$, which assigns
probability $ B^{-\ell}$ to all cylinder sets of length
$\ell$.  If ${\mathbb{Y}}_1,{\mathbb{Y}}_2\subset{\mathbb{Z}}$ are disjoint, and
${\mathbf{ b}}_k\in{\mathcal{ B}}^{{\mathbb{Y}}_k}$, then we define
${\mathbf{ b}}_1\underline{\ }{\mathbf{ b}}_2\in{\mathcal{ B}}^{{\mathbb{Y}}_1\sqcup{\mathbb{Y}}_2}$ by: \
$\left({\mathbf{ b}}_1\underline{\ }{\mathbf{ b}}_2\right)\raisebox{-0.3em}{$\left|_{{\mathbb{Y}}_k}\right.$} = {\mathbf{ b}}_k$, for $k=1,2$.
Thus, ${\left\langle {\mathbf{ b}}_1\underline{\ }{\mathbf{ b}}_2 \right\rangle } = {\left\langle {\mathbf{ b}}_1 \right\rangle }\cap {\left\langle {\mathbf{ b}}_2 \right\rangle }$.
Also, if ${\mathbb{Y}}\subset{\mathbb{Z}}$ and ${\mathbb{X}}={\mathbb{Y}}+{\mathbb{V}}$, then we abuse
notation by letting ${\mathfrak{ G}}:{\mathcal{ B}}^{\mathbb{X}} {{\longrightarrow}} {\mathcal{ B}}^{\mathbb{Y}}$ be the ``local map''
induced by ${\mathfrak{ G}}$.

\subsection{Permutativity and Relative Permutativity}
 A local map ${\mathfrak{ g}}:{\mathcal{ B}}^{{\mathbb{V}}}{{\longrightarrow}}{\mathcal{ B}}$ is {\bf left permutative} if
$ L>0$ and, for every ${\mathbf{ b}} \in {\mathcal{ B}}^{{\left( - L... R \right]}}$ the
map ${\mathcal{ B}}\ni  a\mapsto {\mathfrak{ g}}( a\underline{\ }{\mathbf{ b}}) \in {\mathcal{ B}}$ is a bijection;
\ ${\mathfrak{ g}}$ is {\bf right permutative} if $ R>0$ and, for every ${\mathbf{ b}}
\in {\mathcal{ B}}^{{\left[ - L... R \right)}}$ the map ${\mathcal{ B}}\ni  c\mapsto
{\mathfrak{ g}}({\mathbf{ b}}\underline{\ } c) \in {\mathcal{ B}}$ is a bijection.  ${\mathfrak{ g}}$ is {\bf
permutative} if it is either left- or right-permutative.  A CA on
${\mathcal{ B}}^{\mathbb{Z}}$ is {\bf bipermutative} if it is permutative on both sides;
\ a CA on ${\mathcal{ B}}^{\mathbb{N}}$ is called {\bf bipermutative} if it is
right-permutative.  If $ V= R+ L$ then we 
say ${\mathfrak{ g}}$ is {\bf $ V$-bipermutative}.  A nonhomogeneous CA ${\mathfrak{ G}}$ is
{\bf $ V$-bipermutative} if ${\mathfrak{ g}}_{{\mathsf{ m}}}$ is $ V$-bipermutative for
every ${\mathsf{ m}}\in{\mathbb{M}}$.  A relative CA
${\mathfrak{ F}}:{\mathcal{ A}}^{\mathbb{M}}\times{\mathcal{ C}}^{\mathbb{M}}{{\longrightarrow}}{\mathcal{ A}}^{\mathbb{M}}$ is {\bf
$ V$-bipermutative} if the NHCA ${\mathfrak{ F}}_{\mathbf{ c}}$ is $ V$-bipermutative
for every ${\mathbf{ c}} \in {\mathcal{ C}}^{\mathbb{M}}$.

\medskip

  \medskip         \refstepcounter{thm} {\bf Example \thethm:}  \setcounter{enumi}{\thethm} \begin{list}{(\alph{enumii})}{\usecounter{enumii}} 			{\setlength{\leftmargin}{0em} 			\setlength{\rightmargin}{0em}}   

\item \label{entrXLin}
If ${\mathcal{ B}} = \left({{\mathbb{Z}}_{/n}},+\right)$ and $ V_0<0< V_1$, then a linear
CA with local map \ $\displaystyle {\mathfrak{ g}}({\mathbf{ b}}) \ = \
\sum_{ v= V_0}^{ V_1}  g_v\cdot b_{ v}$ \ 
is left- (resp. right-) permutative iff $ g_{ V_0}$ 
(resp. $ g_{ V_1}$) is relatively prime to $n$.

\item \label{entrx1}
In Example (\ref{relCAx1}), ${\mathfrak{ H}}$ is a $3$-bipermutative CA, while
 ${\mathfrak{ F}}$ is a $3$-bipermutative RCA.

\item \label{entrx2} 
In Example (\ref{relCAx2}),
${\mathfrak{ h}}(c_0,c_1,c_2) = {\mathfrak{ h}}(c_0,c_1) = c_0-c_1$ is right-permutative,
with ${\mathbb{V}}=\{0,1\}$, so ${\mathfrak{ H}}$ is 2-bipermutative as a map on ${\mathcal{ C}}^{\mathbb{N}}$. However, ${\mathfrak{ F}}$ is {\em not} right-permutative.  To see this, write:
\ $\displaystyle
 {\mathfrak{ f}}_{( c_{_{0}}, c_{_{1}}, c_{_{2}})} \left( a_{_{0}}, a_{_{1}}, a_{_{2}}\right)
\ = \ 
  {\mathfrak{ f}}^0_{( c_{_{0}}, c_{_{1}}, c_{_{2}})} ( a_{_{0}})
\ + \ {\mathfrak{ f}}^1_{( c_{_{0}}, c_{_{1}}, c_{_{2}})} ( a_{_{1}})
\ + \ {\mathfrak{ f}}^2_{( c_{_{0}}, c_{_{1}}, c_{_{2}})} ( a_{_{2}})$.
 \[
\mbox{Then:} \hspace{3em}
{\mathfrak{ f}}^2_{( c_{_{0}}, c_{_{1}}, c_{_{2}})} \ = \
  1 +  2^{ c_{2}} +  2^{2 c_{2}} + 2^{3 c_{2}} \  = \
{\left\{ \begin{array}{rcl}                                  4 \pmod{5}&&\mathrm{if}\  c_{_{2}}=0 \ \mathrm{or} \ 4;\\
	 0 \pmod{5}&&\mathrm{if} \  c_{_{2}}=1,2 \ \mathrm{or} \ 3;                                \end{array}  \right.  }
\]
  Thus, ${\mathfrak{ f}}_{( c_{_{0}}, c_{_{1}}, c_{_{2}})}$ is
right-permutative if and only if $ c_{_{2}}=0$ or $4$.
 	\hrulefill\end{list}   			 

The following results extend well-known properties of permutative
cellular automata \cite{HedlundCA} to the nonhomogeneous case; the
proofs are similar, and are left to the reader.

\begin{lemma}\label{permutative.filling}  
Let $ J< K$, and $\ell= K- J$.  Let ${\mathfrak{ G}}$ be an NHCA, and
let ${\mathbf{ d}}\in{\mathcal{ B}}^{\left[  J... K \right)}$.
\begin{enumerate}
\item  If ${\mathfrak{ G}}$ is right-permutative, then
for all $ \ {\mathbf{ b}}\in{\mathcal{ B}}^{\left[  J- L... J+ R \right)}$, there is a
 unique ${\mathbf{ c}} \in {\mathcal{ B}}^{\left[  J+ R... K+ R \right)}$ so that
${\mathfrak{ G}}({\mathbf{ b}}\underline{\ }{\mathbf{ c}}) \ = \ {\mathbf{ d}}$.

\item If ${\mathfrak{ G}}$ is left-permutative, then
for all $ \ {\mathbf{ b}}\in{\mathcal{ B}}^{\left[  K- L... K+ R \right)}$, there is a
 unique ${\mathbf{ a}} \in {\mathcal{ B}}^{\left[  J- L... K- L \right)}$ so that
${\mathfrak{ G}}({\mathbf{ a}}\underline{\ }{\mathbf{ b}}) \ = \ {\mathbf{ d}}$.

\item If ${\mathfrak{ G}}$ is bipermutative, then for any $ j\in{\left[  J... K \right)}$ and
${\mathbf{ b}}\in{\mathcal{ B}}^{\left[  j- L... j+ R \right)}$, there are unique
${\mathbf{ a}} \in {\mathcal{ B}}^{\left[  J- L... j- L \right)}$ and ${\mathbf{ c}} \in
{\mathcal{ B}}^{\left[  j+ R... K+ R \right)}$ so that
${\mathfrak{ G}}({\mathbf{ a}}\underline{\ }{\mathbf{ b}}\underline{\ }{\mathbf{ c}}) \ = \ {\mathbf{ d}}$.\hrulefill\ensuremath{\Box}
\end{enumerate}
 \end{lemma}

\begin{cor}\label{permutative.haar}  
If  ${\mathfrak{ G}}$ is permutative, then $\eta_{_{{\mathcal{ B}}}}{}$ is ${\mathfrak{ G}}$-invariant.\hrulefill\ensuremath{\Box}
 \end{cor}

\subsection{Measurable Entropy}

  Suppose $({\mathbf{ Y}},{\mathcal{ Y}},\mu)$ is a probability space, ${\mathcal{ Q}}$
is a finite set, and ${\mathbf{ Q}}:{\mathbf{ Y}}{{\longrightarrow}}{\mathcal{ Q}}$ is measurable; \
 we say ${\mathbf{ Q}}$ is a {\bf partition} of ${\mathbf{ Y}}$, indexed
by ${\mathcal{ Q}}$.  If $\rho={\mathbf{ Q}}(\mu)\in\ensuremath{{\mathcal{ M}}\left[{\mathcal{ Q}}\right] }$, then let
\ $\displaystyle  h({\mathbf{ Q}};\mu)  =  -\sum_{ q\in{\mathcal{ Q}}} \rho[ q]\cdot \log\left(\rho[ q]\right)$.
If ${\mathbf{ Q}}_k:{\mathbf{ Y}}{{\longrightarrow}}{\mathcal{ Q}}_k$ for $k=1,2$, then let
${\mathbf{ Q}}_1 \vee {\mathbf{ Q}}_2:{\mathbf{ Y}}{{\longrightarrow}} {\mathcal{ Q}}_1 \times {\mathcal{ Q}}_2$ 
be the partition mapping
$ y\in{\mathbf{ Y}}$ to $\left({\mathbf{ Q}}_1( y),{\mathbf{ Q}}_2( y)\right)$.  
If ${\sf G}:{\mathbf{ Y}}\!\longrightarrow{\mathbf{ Y}}\!$ is $\mu$-preserving, then
define ${\sf G}{\mathbf{ Q}} \ = \ {\mathbf{ Q}}\circ {\sf G}$.  If ${\sf G}$ is invertible (respectively noninvertible),
let ${\mathbb{T}}={\mathbb{Z}}$ (respectively ${\mathbb{T}}={\mathbb{N}}$); \ for any ${\mathbb{I}}\subset{\mathbb{T}}$,
define
$\displaystyle {\sf G}^{{\mathbb{I}}}\, {\mathbf{ Q}} \ = \ \bigvee_{i\in{\mathbb{I}}} {\sf G}^i\, {\mathbf{ Q}}:{\mathbf{ Y}}{{\longrightarrow}}{\mathcal{ Q}}^{\mathbb{I}}$. \ 
The {\bf ${\sf G}$-entropy} of ${\mathbf{ Q}}$ is the limit: \ \ 
$\displaystyle
  h({\mathbf{ Q}},{\sf G},\mu) \ = \ 
\lim_{N{\rightarrow}{\infty}}\frac{1}{N}\   h\left({\sf G}^{{\left[ 0...N \right)}} {\mathbf{ Q}};\ \mu \right)$.

  Let $\Sigma\left({\sf G}^{\mathbb{T}}\,{\mathbf{ Q}}\right)$ be the smallest
$\sigma$-algebra for which the function ${\sf G}^{\mathbb{T}}\,{\mathbf{ Q}}$ is
measurable.  ${\mathbf{ Q}}$ is a {\bf generator} for ${\sf G}$ if
$\Sigma\left({\sf G}^{\mathbb{T}}\,{\mathbf{ Q}}\right)$ is $\mu$-{\bf dense} in
${\mathcal{ Y}}$, meaning that, for all ${\mathbf{ V}}\in{\mathcal{ Y}}$, there is
${\mathbf{ W}}\in\Sigma\left({\sf G}^{\mathbb{T}}\,{\mathbf{ Q}}\right)$ so that
$\mu\left({\mathbf{ V}}\Delta{\mathbf{ W}}\right)=0$; \ the (measurable) {\bf entropy} of the
system $({\mathbf{ Y}},{\mathcal{ Y}},\mu;\ {\sf G})$ is then defined:
\ $ h({\mathbf{ Y}},{\mathcal{ Y}},\mu;\ {\sf G}) \ =
\  h({\mathbf{ Q}},{\sf G},\mu)$ \ (independent of the choice of
generator).

 Suppose that $\{{\mathfrak{ R}}_n\}_{n=0}^{\infty}$ is a collection of NHCA on
${\mathcal{ B}}^{\mathbb{M}}$.  For any $n\in{\mathbb{N}}$, define ${\mathfrak{ R}}^{(n)} = {\mathfrak{ R}}_{n-1}
\circ {\mathfrak{ R}}_{n-2}\circ\ldots\circ{\mathfrak{ R}}_0$.  If
${\mathbf{ Q}}$ is a partition, then, for all $N\in{\mathbb{N}}$, define $\displaystyle
{\mathfrak{ R}}^{\left[ 0... N \right)} {\mathbf{ Q}} \ = \ \bigvee_{n=0}^{N-1} {\mathfrak{ R}}^{(n)} {\mathbf{ Q}}$. \
Thus, if ${\mathfrak{ R}}_n={\mathfrak{ G}}$ for all $n\in{\mathbb{N}}$, then ${\mathfrak{ R}}^{(n)} = {\mathfrak{ G}}^{n}$, and
${\mathfrak{ R}}^{\left[ 0... N \right)} {\mathbf{ Q}} \ = \ {\mathfrak{ G}}^{\left[ 0... N \right)} {\mathbf{ Q}}$.

 Say that ${\mathbf{ Q}}_k:{\mathbf{ Y}}{{\longrightarrow}}{\mathcal{ Q}}_k$ ($k=1,2$)
are {\bf equivalent} if each is measurable with respect to the other.
Then, for any $\mu\in\ensuremath{{\mathcal{ M}}\left[{\mathbf{ Y}}\right] }$, \  $ h({\mathbf{ Q}}_1;\mu)\ =\  h({\mathbf{ Q}}_2;\mu)$.

\begin{prop}\label{permutative.entropy}  
Let ${\mathfrak{ R}}_n$ be $ V$-bipermutative, for all $n\in{\mathbb{N}}$.
Let $\mu\in\ensuremath{{\mathcal{ M}}\left[{\mathcal{ B}}^{\mathbb{M}}\right] }$
and let ${\mathbf{ Q}} = {\mathbf{ pr}_{{{\left[ - L... R \right)}}}}:{\mathcal{ B}}^{\mathbb{M}}{{\longrightarrow}}{\mathcal{ B}}^{{\left[ - L... R \right)}}$.
Then:
\begin{enumerate}
 \item   ${\mathfrak{ R}}^{\left[ 0... N \right)} {\mathbf{ Q}}$  and ${\mathbf{ pr}_{{{\left[ -N  L... N  R \right)}}}}$ are equivalent.  Thus,
 $ h\left({\mathfrak{ R}}^{\left[ 0... N \right)} {\mathbf{ Q}};\ \mu\right) =  h\left({\mathbf{ pr}_{{{\left[ -N  L... N  R \right)}}}};\ \mu\right).$

  \item  $\Sigma\left({\mathfrak{ R}}^{\mathbb{N}} {\mathbf{ Q}}\right)$ is the Borel sigma-algebra
of  ${\mathcal{ B}}^{\mathbb{M}}$.

\medskip

\hspace{-3em} In particular, suppose ${\mathfrak{ R}}_n={\mathfrak{ G}}$, for all $n\in{\mathbb{N}}$. 
Then:
  \item ${\mathbf{ Q}}$ is a $({\mathfrak{ G}},\mu)$-generator.

  \item If $\mu$ is $ {{{\boldsymbol{\sigma}}}^{}} $-invariant and ${\mathfrak{ G}}$-invariant, then $ h\left({\mathfrak{ G}};\ \mu\right) \ = \  V\cdot  h\left( {{{\boldsymbol{\sigma}}}^{}} ;\ \mu\right)$.

 In particular, $ h\left({\mathfrak{ G}};\ \eta_{_{{\mathcal{ B}}}}{}\right) \ = \  V\cdot \log( B)$.
\end{enumerate}  
 \end{prop}
\bprf  {\bf Part 1} is proved by repeated application of {\bf Part 3}
from Lemma \ref{permutative.filling}.  The other statements then follow.
{\tt \hrulefill $\Box$ }\end{list}  \medskip  
{\bf Remark:} By combining {\bf Part 4} with Example (\ref{entrXLin}), we
recover the previously computed \cite{DamicoManziniMargara} entropy of
linear CA on $\left(\left({{\mathbb{Z}}_{/p}}\right)^{\mathbb{Z}},\eta_{_{{\mathcal{ B}}}}{}\right)$.

\subsection{Relative Entropy
\label{S:Rel.Entropy}}

  Let ${\mathbf{ Y}}={\mathbf{ X}}\times{\mathbf{ Z}}$, \
${\sf G}={\sf F}\star{\sf H}$, and $\lambda\in\ensuremath{{\mathcal{ M}}\left[{\mathbf{ X}}\right] }$ be
as in \S\ref{S:rel.dyn.sys}.  If $\nu\in\ensuremath{{\mathcal{ M}}\left[{\mathbf{ Z}}\right] }$ is
${\sf H}$-invariant, then $\mu=\lambda\otimes\nu$ is
${\sf G}$-invariant.  For any $ z\in{\mathbf{ Z}}$ and $n\in{\mathbb{N}}$,
define \ $\displaystyle {\sf F}^{(n)}_ z \ = \
{\sf F}_{{\sf H}^{n-1}( z)} \circ\ldots \circ
{\sf F}_{{\sf H}^2( z)} \circ {\sf F}_{{\sf H}( z)} \circ
{\sf F}_ z$. If ${\mathbf{ P}}:{\mathbf{ X}}{{\longrightarrow}}{\mathcal{ P}}$ is a partition,
then for all $ z \in {\mathbf{ Z}}$, let \ ${\sf G}_ z^{{\left[ 0... N \right)}}
{\mathbf{ P}}:{\mathbf{ X}}{{\longrightarrow}}{\mathcal{ P}}^{N}$ be the partition: $ \displaystyle
{\sf G}_ z^{{\left[ 0... N \right)}} {\mathbf{ P}} \ = \ \bigvee_{n=0}^{N-1} {\sf F}^{(n)}_ z {\mathbf{ P}}$, and define
\ $\displaystyle   h\left({\mathbf{ P}}; \ {\sf F},\lambda/ \ {\sf H}, \nu \right)
 \ = \
\lim_{N{\rightarrow}{\infty}}\frac{1}{N}\int_{{\mathbf{ Z}}}
 h\left({\sf G}_ z^{{\left[ 0... N \right)}} {\mathbf{ P}}, \ \lambda \right) \ d\nu[ z]$.
The {\bf relative entropy} of  ${\sf F}$ over ${\sf H}$
is defined:  $\displaystyle
  h\left({\sf F},\lambda / {\sf H},\nu\right)
\ = \ \sup_{{\mathbf{ P}}} \  h\left({\mathbf{ P}}; \ {\sf F},\lambda/ \ {\sf H}, \nu\right)$, where 
the supremum is taken over all measurable partitions ${\mathbf{ P}}$ of
 ${\mathbf{ X}}$.
\begin{thm}\label{abramov.rokhlin}  
$\displaystyle  h\left({\sf G}, \ \lambda \times \nu\right) \ \  =  \ \  
 h\left({\sf H},\nu\right) \  + \   h\left({\sf F},\lambda / {\sf H},\nu\right)$.
\hrulefill {\rm(Abramov \& Rokhlin) \cite{AbramovRokhlin,Petersen}}.
 \end{thm}

\begin{thm}\label{entropy.skewprod}  
 Let ${\mathcal{ B}}={\mathcal{ A}}\times{\mathcal{ C}}$. 
 Let $\lambda\in\ensuremath{{\mathcal{ M}}\left[{\mathcal{ A}}^{\mathbb{M}}\right] }$, $\nu\in\ensuremath{{\mathcal{ M}}\left[{\mathcal{ C}}^{\mathbb{M}}\right] }$,
and  $\mu = \lambda\otimes\nu\in\ensuremath{{\mathcal{ M}}\left[{\mathcal{ B}}^{\mathbb{M}}\right] }$ be
$ {{{\boldsymbol{\sigma}}}^{}} $-invariant.  Suppose ${\mathfrak{ G}} =
{\mathfrak{ F}} \,\star\, {\mathfrak{ H}}$, where
 ${\mathfrak{ F}}:{\mathcal{ A}}^{\mathbb{M}}\times{\mathcal{ C}}^{\mathbb{M}}{{\longrightarrow}} {\mathcal{ A}}^{\mathbb{M}}$ is
$\lambda$-preserving and $ V$-bipermutative, while
${\mathfrak{ H}}:{\mathcal{ C}}^{\mathbb{M}}\!\longrightarrow{\mathcal{ C}}^{\mathbb{M}}\!$ is $\nu$-preserving and
$ W$-bipermutative.  Then:\\
{\bf(1)} $ h\left({\mathfrak{ F}},\lambda/\ {\mathfrak{ H}},\nu\right)
 \ = \   V \cdot  h\left(\lambda, {{{\boldsymbol{\sigma}}}^{}} \right)$;\\
{\bf(2)} $ h\left({\mathfrak{ H}},\nu\right) \ = \ 
 W \cdot  h\left(\nu, {{{\boldsymbol{\sigma}}}^{}} \right)$; \\
{\bf(3)} $ h({\mathfrak{ G}},\mu) \ = \  V \cdot  h\left(\lambda, {{{\boldsymbol{\sigma}}}^{}} \right) +  W \cdot  h\left(\nu, {{{\boldsymbol{\sigma}}}^{}} \right)$.

 \end{thm}
\bprf To see {\bf(1)}, \
let ${\mathbf{ Q}} =
{\mathbf{ pr}_{{{\left[ - L... R \right)}}}}:{\mathcal{ A}}^{\mathbb{M}}{{\longrightarrow}}{\mathcal{ A}}^{{\left[ - L... R \right)}}$.
 Fix ${\mathbf{ c}} \in {\mathcal{ C}}^{\mathbb{M}}$, and, for all $n\in{\mathbb{N}}$, let 
${\mathfrak{ R}}_n = {\mathfrak{ F}}_{{\mathfrak{ H}}^n({\mathbf{ c}})}$.  Then ${\mathfrak{ R}}^{(n)} = {\mathfrak{ F}}^{(n)}_{\mathbf{ c}}$,
and ${\mathfrak{ R}}^{{\left[ 0... N \right)}} {\mathbf{ Q}} = {\mathfrak{ G}}^{{\left[ 0... N \right)}}_{\mathbf{ c}} {\mathbf{ Q}}$, so
{\bf Part 1}
of Proposition \ref{permutative.entropy} says: 
\ $\displaystyle  h\left({\mathfrak{ G}}_{\mathbf{ c}}^{{\left[ 0... N \right)}} {\mathbf{ Q}}; \ \lambda \right)
\ = \ 
 h\left({\mathbf{ pr}_{{{\left[ -N  L... N  R \right)}}}};\ \lambda\right)$.

  Let ${\mathbf{ P}}:{\mathcal{ A}}^{\mathbb{M}}{{\longrightarrow}}{\mathcal{ P}}$ be any other partition of
${\mathcal{ A}}^{\mathbb{M}}$, and fix $\epsilon>0$.  Then, by {\bf Part 2} of Proposition
\ref{permutative.entropy}, there is some $M=M({\mathbf{ c}})>0$ so that, for
all $N\in{\mathbb{N}}$, \ $ h\left({\mathfrak{ G}}_{\mathbf{ c}}^{{\left[ 0... N \right)}} {\mathbf{ P}}; \
\lambda \right)
\ < \ N\epsilon \ + \ 
 h\left({\mathfrak{ G}}_{\mathbf{ c}}^{{\left[ 0... N+M \right)}} {\mathbf{ Q}}; \ \lambda \right)$. 
If ${{\sf card}\left[{\mathcal{ P}}\right]}=P$, then also,
$ h\left({\mathfrak{ G}}_{\mathbf{ c}}^{{\left[ 0... N \right)}} {\mathbf{ P}}; \ \lambda \right) 
< N\log(P)$.
Find $M$ so that $\mu[{\mathcal{ D}}]<\epsilon/\log(P)$, where
${\mathcal{ D}}={\left\{ {\mathbf{ c}}\in{\mathcal{ C}}^{\mathbb{M}} \; ; \; M({\mathbf{ c}})>M \right\} }$.
Then
\begin{eqnarray*}
 \int_{{\mathcal{ C}}^{\mathbb{M}}}
 h\left({\mathfrak{ G}}_{\mathbf{ c}}^{{\left[ 0... N \right)}} {\mathbf{ P}}; \ \lambda \right) \ d\nu[{\mathbf{ c}}]
& = & 
 \int_{{\mathcal{ D}}}
 h\left({\mathfrak{ G}}_{\mathbf{ c}}^{{\left[ 0... N \right)}} {\mathbf{ P}}; \ \lambda \right) \ d\nu[{\mathbf{ c}}]
\ + \ 
 \int_{{\mathcal{ C}}^{\mathbb{M}}\setminus {\mathcal{ D}}}
 h\left({\mathfrak{ G}}_{\mathbf{ c}}^{{\left[ 0... N \right)}} {\mathbf{ P}}; \ \lambda \right) \ d\nu[{\mathbf{ c}}]
\\
&<&
\left(N\log(P)\frac{\epsilon}{\log(P)} \right) \ \ + \  N\epsilon \ + \ 
 h\left(\rule[-0.5em]{0em}{1em} \mathbf{pr}_{\left[-(N+M)  L...(N+M)R \right)};\ \lambda\right).
\end{eqnarray*}
Thus,
\ \ $\displaystyle  h \left(\mathbf{P}; \ \mathfrak{F}, \lambda\ / \ \mathfrak{H},\nu  \right)  
 < \ \    
2 \epsilon + 
\lim_{N{\rightarrow}{\infty}}\frac{1}{N} h\left(\mathbf{pr}_{\left[-N  L...N R\right)};\ \lambda\right)
 = \ \   2 \epsilon \ + \ V\cdot h\left(\boldsymbol{\sigma};\ \lambda\right)$.
 Take the supremum
over all $\mathbf{P}$ to conclude: \ 
$\displaystyle h\left(\mathfrak{F}, \lambda\ / \ \mathfrak{H},\nu\right) \;  < \, 
2\epsilon \, + \,   V\cdot h\left(\boldsymbol{\sigma};\ \lambda\right)$.  Now let $\epsilon\rightarrow0$.

{\bf(2)} follows from {\bf Part 4} of
Proposition \ref{permutative.entropy} and {\bf(3)} follows from Theorem \ref{abramov.rokhlin}.
{\tt \hrulefill $\Box$ }\end{list}  \medskip  

  Note that, in Theorem \ref{entropy.skewprod}, ${\mathcal{ B}}$ need not be a
group, nor ${\mathfrak{ G}}$ a multiplicative cellular automaton.  However,
Theorem \ref{MCA.structure} provides a natural skew product decomposition
in this case.
 
  \medskip         \refstepcounter{thm}                     \begin{list}{} 			{\setlength{\leftmargin}{1em} 			\setlength{\rightmargin}{0em}}                     \item {\bf Example \thethm:} In Example (\ref{entrx1}), $ V= W=2$,
${{\sf card}\left[{\mathcal{ A}}\right]}=5$ and ${{\sf card}\left[{\mathcal{ C}}\right]}=4$; \ thus $
 h({\mathfrak{ G}},\eta_{_{{\mathcal{ B}}}}{}) = 2\cdot \log_2(5)+2\cdot \log_2(4) = 2\log_2(5)
+ 4$.    	\hrulefill\end{list}

\section{Convergence of Measures
\label{S:measure}}

  Endow $\ensuremath{{\mathcal{ M}}\left[{\mathcal{ B}}^{\mathbb{M}}\right] }$ with the weak* topology induced by ${\mathbf{ C}}\left({\mathcal{ B}}^{\mathbb{M}}; \
{\mathbb{C}}\right)$, the space of continuous, complex-valued functions.
The uniformly distributed Bernoulli measure $\eta_{_{{\mathcal{ B}}}}{} \in
\ensuremath{{\mathcal{ M}}\left[{\mathcal{ B}}^{\mathbb{M}}\right] }$ is the Haar measure on ${\mathcal{ B}}^{\mathbb{M}}$ as a compact group,
and is invariant under the action of any left- or right-permutative
MCA (Lemma \ref{permutative.haar}).  Thus, if
$\mu\in\ensuremath{{\mathcal{ M}}\left[{\mathcal{ B}}^{\mathbb{M}}\right] }$ is some initial measure, then
$\eta_{_{{\mathcal{ B}}}}{}$ is a natural candidate for the (Ces\`aro) limit of ${\mathfrak{ G}}^n \mu$ 
as $n{\rightarrow}{\infty}$.  Since $\eta_{_{{\mathcal{ B}}}}{}$ is the measure of maximal
entropy on ${\mathcal{ B}}^{\mathbb{M}}$, such limiting behaviour is a sort of
``asymptotic randomization'' of ${\mathcal{ B}}^{\mathbb{M}}$.  When ${\mathcal{ B}}$ is abelian, and
${\mathfrak{ G}}$ is an affine CA, the
Ces\`aro  convergence of measures to $\eta_{_{{\mathcal{ B}}}}{}$ is somewhat understood
\cite{MaassMartinezII,Lind,PivatoYassawi1,PivatoYassawi2,FerMaassMartNey}; \
we now extend these results to nonabelian MCA.

\subparagraph*{Harmonic Mixing and Diffusion:}
  The {\bf characters} of an abelian group $({\mathcal{ A}},+)$ are the
continuous homomorphisms from ${\mathcal{ A}}$ into the unit circle group
$\left({{{\mathbb{T}}}^{1}},\cdot\right)\subset{\mathbb{C}}$.  The set of characters forms a
group, denoted ${\widehat{{\mathcal{ A}}}}$.  For example, if ${\mathcal{ A}}={{\mathbb{Z}}_{/n}}$, then
every $\chi\in{\widehat{{\mathcal{ A}}}}$ has the form $\chi(a)
= \exp\left(\frac{2\pi{\mathbf{ i}}}{n} c\cdot a\right)$, where
$c\in{{\mathbb{Z}}_{/n}}$ is a constant coefficient, and the product $c\cdot a$ is
computed mod $n$.

 For any ${\boldsymbol{\chi }}\in {\widehat{{\mathcal{ A}}^{\mathbb{M}}}}$ there is some finite ${\mathbb{K}}\subset{\mathbb{M}}$
and, for each ${\mathsf{ k}}\in{\mathbb{K}}$, a nontrivial $\chi_{\mathsf{ k}}\in{\widehat{{\mathcal{ A}}}}$,
so that, if ${\mathbf{ a}}\in{\mathcal{ A}}^{\mathbb{M}}$, then $\displaystyle {\boldsymbol{\chi }}({\mathbf{ a}})
\ = \ \prod_{{\mathsf{ k}}\in{\mathbb{K}}} \chi_{\mathsf{ k}} (a_{\mathsf{ k}})$; \ we indicate this:
``$\displaystyle{\boldsymbol{\chi }}=\bigotimes_{{\mathsf{ k}}\in{\mathbb{K}}} \chi_{\mathsf{ k}}$''.  For example,
if ${\mathcal{ A}}={{\mathbb{Z}}_{/n}}$, this means there is a collection of nonzero
{coefficients} ${\left[c_{\mathsf{ k}}  |_{{\mathsf{ k}}\in{\mathbb{K}}}^{} \right]}$ so that if
${\mathbf{ a}}\in{\mathcal{ A}}^{\mathbb{M}}$, then $\displaystyle{\boldsymbol{\chi }}({\mathbf{ a}})
\ = \ \exp\left(\frac{2\pi{\mathbf{ i}}}{n} \sum_{{\mathsf{ k}}\in{\mathbb{K}}} c_{\mathsf{ k}}\cdot a_{\mathsf{ k}}\right)$.

The {\bf rank} of ${\boldsymbol{\chi }}\in {\widehat{{\mathcal{ A}}^{\mathbb{M}}}}$ is the cardinality of ${\mathbb{K}}$.  Let $\mu
\in\ensuremath{{\mathcal{ M}}\left[{\mathcal{ A}}^{\mathbb{M}}\right] }$; we will use the notation
${\left\langle {\boldsymbol{\chi }},\mu \right\rangle } \ = \ \int_{{\mathcal{ A}}^{\mathbb{M}}} {\boldsymbol{\chi }} \ d\mu$.
We call $\mu$ {\bf harmonically mixing}
(and write ``$\,\mu\in{\mathcal{ H}}{\mathcal{ M}}\left[{\mathcal{ A}}^{\mathbb{M}}\right]$'')
 if, for all 
$\epsilon>0$,\ there is  $r\in{\mathbb{N}}$ so that, if
${\boldsymbol{\chi }}\in{\widehat{{\mathcal{ A}}^{\mathbb{M}}}}$ and ${{\sf rank}\left[{\boldsymbol{\chi }}\right]}>r$, then
$\left|{\left\langle \mu,\chi \right\rangle }\right|<\epsilon$.  For example, most Bernoulli measures
\cite{PivatoYassawi1} and Markov random fields
\cite{PivatoYassawi2} on ${\mathcal{ A}}^{\left({\mathbb{Z}}^D\right)}$ are harmonically mixing.

If ${\mathfrak{ G}}\in{{{\mathbf{ E}}{\mathbf{ n}}{\mathbf{ d}}_{} \left[{\mathcal{ A}}^{\mathbb{M}}\right]}}$, then for any ${\boldsymbol{\chi }}\in{\widehat{{\mathcal{ A}}^{\mathbb{M}}}}$, the
map ${\boldsymbol{\chi }}\circ {\mathfrak{ G}}$ is also a character.  If ${\mathbb{J}}\subset{\mathbb{N}}$, then
${\mathfrak{ G}}$ is {\bf ${\mathbb{J}}$-diffusive} if, for every ${\boldsymbol{\chi }}\in{\widehat{{\mathcal{ A}}^{\mathbb{M}}}}$, \
$\displaystyle\lim_{{j{\rightarrow}{\infty}}\atop{j\in{\mathbb{J}}}}
{{\sf rank}\left[{\boldsymbol{\chi }}\circ {\mathfrak{ G}}^j\right]} \ = \ {\infty}$.  If ${\mathbb{J}}={\mathbb{Z}}$, then we just say
${\mathfrak{ G}}$ is {\bf diffusive}; if ${\mathbb{J}}\subset{\mathbb{N}}$ is a subset of Ces\`aro 
density 1, then ${\mathfrak{ G}}$ is {\bf diffusive in density}.

  \medskip         \refstepcounter{thm}                     \begin{list}{} 			{\setlength{\leftmargin}{1em} 			\setlength{\rightmargin}{0em}}                     \item {\bf Example \thethm:}  \label{X:LCA.diffuse}
Let $({\mathcal{ A}},+)$ be a finite abelian group,  \ ${\mathbb{M}}={\mathbb{Z}}^D$,
and let ${\mathfrak{ L}}$ be a linear CA with local map 
$\displaystyle {\mathfrak{ l}}({\mathbf{ a}}) = \sum_{{\mathsf{ v}}\in{\mathbb{V}}} \ell_{\mathsf{ v}} a_{\mathsf{ v}}$, where
$\ell_{\mathsf{ v}}\in{\mathbb{Z}}$ is relatively prime to ${{\sf card}\left[{\mathcal{ A}}\right]}$ for all ${\mathsf{ v}}\in{\mathbb{V}}$.
Then ${\mathfrak{ L}}$ is diffusive in density \cite{PivatoYassawi2}.  	\hrulefill\end{list}   			  
 If ${\mathfrak{ G}}$ is
${\mathbb{J}}$-diffusive and $\mu\in{\mathcal{ H}}{\mathcal{ M}}\left[{\mathcal{ A}}^{\mathbb{M}}\right]$, then 
Theorem 12 of \cite{PivatoYassawi1} says $\displaystyle
\mathrm{wk^*}\!\!\lim_{{\mathbb{J}}\ni j{\rightarrow}{\infty}} {\mathfrak{ G}}^n \mu \ = \ \eta_{_{{\mathcal{ B}}}}{}$. 
In particular, if ${\mathfrak{ G}}$ is diffusive in density, then
 the Ces\`aro  average  weak*-converges to Haar measure:
\begin{equation}
\label{eqn:cesaro.to.haar}
\mathrm{wk^*}\!\!\lim_{N{\rightarrow}{\infty}}\frac{1}{N}
\sum_{n=1}^N {\mathfrak{ G}}^n \mu \ \ = \ \ \eta_{_{{\mathcal{ B}}}}{}.
\end{equation}
  To extend these results to multiplicative cellular automata, we
need a version of diffusion applicable to affine relative 
cellular automata. An {\bf affine endomorphism} of ${\mathcal{ A}}^{\mathbb{M}}$ is a self-map
${\mathfrak{ G}}:{\mathcal{ A}}^{\mathbb{M}}\!\longrightarrow{\mathcal{ A}}^{\mathbb{M}}\!$ of the form ${\mathfrak{ G}}({\mathbf{ a}}) = {\mathbf{ c}} + {\mathfrak{ L}}({\mathbf{ a}})$, where
${\mathfrak{ L}}\in{{{\mathbf{ E}}{\mathbf{ n}}{\mathbf{ d}}_{} \left[{\mathcal{ A}}^{\mathbb{M}}\right]}}$
and ${\mathbf{ c}}\in{\mathcal{ A}}^{\mathbb{M}}$ is constant.   The set
of affine endomorphisms is denoted ${\widetilde{\mathbf{End}}\left[{\mathcal{ A}}^{\mathbb{M}}\right]}$.

An {\bf affine character} is a function
${\boldsymbol{\alpha}}:{\mathcal{ A}}^{\mathbb{M}}{{\longrightarrow}}{{{\mathbb{T}}}^{1}}$ of the form ${\boldsymbol{\alpha}}({\mathbf{ a}}) = c\cdot
{\boldsymbol{\chi }}({\mathbf{ a}})$, where $c\in{{{\mathbb{T}}}^{1}}$ is a constant, and
${\boldsymbol{\chi }}\in{\widehat{{\mathcal{ A}}^{\mathbb{M}}}}$.  The set of affine characters of
${\mathcal{ A}}^{\mathbb{M}}$ is denoted ${\widetilde{{\mathcal{ A}}^{\mathbb{M}}}}$.  For example, if
${\boldsymbol{\zeta}}\in{\widehat{{\mathcal{ A}}^{\mathbb{M}}}}$, and ${\mathfrak{ G}}\in{\widetilde{\mathbf{End}}\left[{\mathcal{ A}}^{\mathbb{M}}\right]}$, then
${\boldsymbol{\zeta}}\circ{\mathfrak{ G}} \in
{\widetilde{{\mathcal{ A}}^{\mathbb{M}}}}$.
The {\bf rank} of ${\boldsymbol{\alpha}}=c\cdot{\boldsymbol{\chi }}$ is the rank of ${\boldsymbol{\chi }}$.  If $\mu
\in\ensuremath{{\mathcal{ M}}\left[{\mathcal{ A}}^{\mathbb{M}}\right] }$ is harmonically mixing, with $\epsilon$ and $r$ as before,
then it follows that $\left|{\left\langle \mu,{\boldsymbol{\alpha}} \right\rangle }\right|<\epsilon$ for any
${\boldsymbol{\alpha}}\in{\widetilde{{\mathcal{ A}}^{\mathbb{M}}}}$ with ${{\sf rank}\left[{\boldsymbol{\alpha}}\right]}>r$.

\subparagraph*{Relative Diffusion:}
Suppose ${\mathcal{ B}}={\mathcal{ A}}\star{\mathcal{ C}}$, where
${\mathcal{ A}}$ is abelian, and let ${\mathfrak{ G}}={\mathfrak{ F}}\star{\mathfrak{ H}}$ as in
Theorem \ref{MCA.structure}.  For any
${\mathbf{ c}}\in{\mathcal{ C}}^{\mathbb{M}}$,
the fibre map ${\mathfrak{ F}}_{\mathbf{ c}}$ is an affine endomorphism; \ 
we say that ${\mathfrak{ F}}$ is an {\em affine relative cellular automaton} (ARCA).
For any $j\in{\mathbb{N}}$, \  ${\mathfrak{ G}}^j =
{\mathfrak{ F}}^{(j)} \star {\mathfrak{ H}}^j$, where ${\mathfrak{ F}}^{(j)}$ is
another ARCA, so ${\boldsymbol{\alpha}}\circ{\mathfrak{ F}}^{(j)}_{\mathbf{ c}}
\in{\widetilde{{\mathcal{ A}}^{\mathbb{M}}}}$ for any ${\mathbf{ c}}\in{\mathcal{ C}}^{\mathbb{M}}$ and
${\boldsymbol{\alpha}}\in{\widetilde{{\mathcal{ A}}^{\mathbb{M}}}}$.  We say ${\mathfrak{ G}}$ is {\bf relatively
${\mathbb{J}}$-diffusive} if $\displaystyle\lim_{{\mathbb{J}}\ni j{\rightarrow}{\infty}}
{{\sf rank}\left[{\boldsymbol{\alpha}}\circ{\mathfrak{ F}}^{(j)}_{\mathbf{ c}}\right]} \ = \ {\infty}$ \
for every ${\mathbf{ c}}\in{\mathcal{ C}}^{\mathbb{M}}$ and ${\boldsymbol{\alpha}}\in{\widetilde{{\mathcal{ A}}^{\mathbb{M}}}}$.
If $\nu\in\ensuremath{{\mathcal{ M}}\left[{\mathcal{ C}}^{\mathbb{M}}\right] }$, and
${\mathbb{J}}\subset{\mathbb{N}}$, then ${\mathfrak{ G}}$ is {\bf $\nu$-relatively
${\mathbb{J}}$-diffusive} if, 
\begin{equation}
\label{eqn:mu.rel.dif}
\forall {\boldsymbol{\alpha}}\in{\widetilde{{\mathcal{ A}}^{\mathbb{M}}}}, \ 
\forall r>0, \ \ \ 
  \lim_{{\mathbb{J}}\ni j {\rightarrow}{\infty}} \nu{\left\{ {\mathbf{ c}}\in{\mathcal{ C}}^{\mathbb{M}} \; ; \; {{\sf rank}\left[{\boldsymbol{\alpha}}\circ{\mathfrak{ F}}^{(j)}_{\mathbf{ c}}\right]} \leq r \right\} } \ = \ 0.
 \end{equation}
Clearly, relative diffusion implies $\nu$-relative diffusion for any
$\nu\in\ensuremath{{\mathcal{ M}}\left[{\mathcal{ C}}^{\mathbb{M}}\right] }$.

\begin{prop}\label{X:aff.diffuse}  
 If ${\mathcal{ A}}\subset Z({\mathcal{ B}})$ as in {\bf Part 3} of Proposition
\ref{MCA.structure.prop}, then ${\mathfrak{ F}}$ is relatively ${\mathbb{J}}$-diffusive if and
only if ${\mathfrak{ L}}$ is ${\mathbb{J}}$-diffusive as a linear cellular automaton.
 \end{prop}
\bprf 
For any $N\in{\mathbb{N}}$, define ${\mathfrak{ P}}^{(N)} = 
\displaystyle \sum_{n=0}^{N-1} {\mathfrak{ L}}^n\circ{\mathfrak{ P}}\circ{\mathfrak{ H}}^{N-n-1}$.
If $j\in{\mathbb{J}}$, then ${\mathfrak{ F}}_{\mathbf{ c}}^{(j)} \ = \ {\mathfrak{ L}}^j +
{\mathfrak{ P}}^{(j)}({\mathbf{ c}})$. Thus, for any
${\boldsymbol{\alpha}}\in{\widetilde{{\mathcal{ A}}^{\mathbb{M}}}}$, ${{\sf rank}\left[{\boldsymbol{\alpha}}\circ{\mathfrak{ F}}_{\mathbf{ c}}^j\right]}
\ = \ {{\sf rank}\left[{\boldsymbol{\alpha}}\circ {\mathfrak{ L}}^j\right]}$.
{\tt \hrulefill $\Box$ }\end{list}  \medskip  

\begin{prop}\label{P:Reldiffuse.into.Haar}  
 Let $\lambda\in{\mathcal{ H}}{\mathcal{ M}}\left[{\mathcal{ A}}^{\mathbb{M}}\right]$, $\nu\in\ensuremath{{\mathcal{ M}}\left[{\mathcal{ C}}^{\mathbb{M}}\right] }$,
 and  $\mu = \lambda\otimes\nu \in \ensuremath{{\mathcal{ M}}\left[{\mathcal{ B}}^{\mathbb{M}}\right] }$.  Let $\displaystyle
 {\overline{\nu }} = \mathrm{wk^*}\!\!\lim_{{\mathbb{J}}\ni j{\rightarrow}{\infty}} {\mathfrak{ H}}^j \nu$, \ and let
 $\eta_{_{{\mathcal{ A}}}}$ be the Haar measure on ${\mathcal{ A}}^{\mathbb{M}}$.  If 
${\mathfrak{ G}}$ is $\nu$-relatively ${\mathbb{J}}$-diffusive,
 then $\displaystyle \mathrm{wk^*}\!\!\lim_{{\mathbb{J}}\ni j{\rightarrow}{\infty}} {\mathfrak{ G}}^j \mu \ = \ \eta_{_{{\mathcal{ A}}}}
 \otimes {\overline{\nu }}$.  
 \end{prop}
\bprf We want
  $\displaystyle \lim_{{\mathbb{J}}\ni j{\rightarrow}{\infty}} {\left\langle {\boldsymbol{\beta}}, \ {\mathfrak{ G}}^j \mu  \right\rangle } \ = \
 {\left\langle {\boldsymbol{\beta}}, \ \eta_{_{{\mathcal{ A}}}} \otimes  {\overline{\nu }} \right\rangle }$,  for every
${\boldsymbol{\beta}}\in{\mathbf{ C}}\left({\mathcal{ B}}^{\mathbb{M}}; \ {\mathbb{C}}\right)$.
It suffices to assume ${\boldsymbol{\beta}} = {\boldsymbol{\alpha}}\otimes {\boldsymbol{\phi }}$, where
${\boldsymbol{\alpha}}\in{\mathbf{ C}}\left({\mathcal{ A}}^{\mathbb{M}}; \ {\mathbb{C}}\right)$
 and ${\boldsymbol{\phi }}\in{\mathbf{ C}}\left({\mathcal{ C}}^{\mathbb{M}}; \ {\mathbb{C}}\right)$.  Since
the characters of ${\mathcal{ A}}^{\mathbb{M}}$ form a basis for the Banach
space ${\mathbf{ C}}\left({\mathcal{ A}}^{\mathbb{M}}; \ {\mathbb{C}}\right)$, it suffices to assume
${\boldsymbol{\alpha}}\in {\widehat{{\mathcal{ A}}^{\mathbb{M}}}}$, and  ${\left\| {\boldsymbol{\phi }} \right\|_{{{\infty}}} }   =1$. \
\[
\mbox{Then:} \ \ 
{\left\langle {\boldsymbol{\beta}}, \ \eta_{_{{\mathcal{ A}}}} \otimes  {\overline{\nu }} \right\rangle } = {\left\langle {\boldsymbol{\alpha}},\eta_{_{{\mathcal{ A}}}} \right\rangle }
\cdot {\left\langle {\boldsymbol{\phi }},{\overline{\nu }} \right\rangle }
\ = \ {\left\{ \begin{array}{rcl}                                  {\left\langle {\boldsymbol{\phi }},{\overline{\nu }} \right\rangle } &&\mbox{if ${\boldsymbol{\alpha}}={{{{\mathsf{ 1\!\!1}}}_{{}}}}$}\\
 		0  &&\mbox{if ${\boldsymbol{\alpha}}\neq{{{{\mathsf{ 1\!\!1}}}_{{}}}}$}                                \end{array}  \right.  }
\]
Now,  for all \ ${\mathbf{ a}}\star{\mathbf{ c}}\in{\mathcal{ B}}^{\mathbb{M}}$, \ \ \
${\boldsymbol{\beta}}\circ {\mathfrak{ G}}^j({\mathbf{ a}}\star{\mathbf{ c}}) \ = \ \left( {\boldsymbol{\alpha}}\otimes {\boldsymbol{\phi }} \right)
\left( {\mathfrak{ F}}^{(j)}_{\mathbf{ c}}\left({\mathbf{ a}}\right) \star  {\mathfrak{ H}}^j \left({\mathbf{ c}}\right) \right)
 \ = \ $ $\displaystyle   \left({\boldsymbol{\alpha}} \circ {\mathfrak{ F}}^{(j)}_{\mathbf{ c}}\left({\mathbf{ a}}\right)\right)
 \cdot \left( {\boldsymbol{\phi }} \circ {\mathfrak{ H}}^j \left({\mathbf{ c}}\right)\right)$.  \ \  Thus, \  
$\displaystyle {\left\langle {\boldsymbol{\beta}}, \ {\mathfrak{ G}}^j \mu \right\rangle } \ = \ 
{\left\langle {\boldsymbol{\beta}}\circ {\mathfrak{ G}}^j, \ \mu \right\rangle }
\ = \ $ $\displaystyle \int_{{\mathcal{ C}}^{\mathbb{M}}}
 \left(\rule[-0.5em]{0em}{1em} {\boldsymbol{\phi }} \circ {\mathfrak{ H}}^j \left({\mathbf{ c}}\right) \right)\cdot
{\left\langle {\boldsymbol{\alpha}} \circ {\mathfrak{ F}}^{(j)}_{\mathbf{ c}}, \ \lambda \right\rangle }  \ d\nu\left[{\mathbf{ c}}\right]$.

  If ${\boldsymbol{\alpha}}={{{{\mathsf{ 1\!\!1}}}_{{}}}}$, then this integral is just equal to $\displaystyle
\int_{{\mathcal{ C}}^{\mathbb{M}}} {\boldsymbol{\phi }} \circ {\mathfrak{ H}}^j \left({\mathbf{ c}}\right) \ d\nu\left[{\mathbf{ c}}\right]$, which converges to ${\left\langle {\boldsymbol{\phi }},\ {\overline{\nu }} \right\rangle }$ by hypothesis.
Hence, assume ${\boldsymbol{\alpha}}\neq{{{{\mathsf{ 1\!\!1}}}_{{}}}}$; \ we then want to show that
$\displaystyle\lim_{{\mathbb{J}}\ni j{\rightarrow}{\infty}} {\left\langle {\boldsymbol{\beta}}, \ {\mathfrak{ G}}^j \mu \right\rangle } = 0$.

  Fix $\epsilon>0$.  Since $\lambda\in{\mathcal{ H}}{\mathcal{ M}}\left[{\mathcal{ A}}^{\mathbb{M}}\right]$, find $r>0$ so
that, if ${\boldsymbol{\alpha}}\in{\widetilde{{\mathcal{ A}}^{\mathbb{M}}}}$ and ${{\sf rank}\left[{\boldsymbol{\alpha}}\right]}>r$, then
$\left|{\left\langle {\boldsymbol{\alpha}},\lambda \right\rangle }\right|<\frac{\epsilon}{2}$. \
Let \ 
$\displaystyle {\mathbf{ D}}_j = {\left\{ {\mathbf{ c}}\in{\mathcal{ C}}^{\mathbb{M}} \; ; \; {{\sf rank}\left[{\boldsymbol{\alpha}}\circ{\mathfrak{ F}}^{(j)}_{\mathbf{ c}}\right]} > r \right\} }$,
for every $j\in{\mathbb{J}}$.
By equation (\ref{eqn:mu.rel.dif}), find
$J\in{\mathbb{N}}$ so that,
$\forall j\in{\mathbb{J}}$ with $j>J$, \ 
$\displaystyle \nu\left[{\mathbf{ D}}_j\right] \ > \ 1-\frac{\epsilon}{2}$.  \ Then \ 
$\displaystyle \left|{\left\langle {\boldsymbol{\beta}}, \ {\mathfrak{ G}}^j \mu \right\rangle }\right| \ = \ 
\left|\int_{{\mathcal{ C}}^{\mathbb{M}}}
  {\boldsymbol{\phi }} \circ {\mathfrak{ H}}^j \left({\mathbf{ c}}\right) \cdot
{\left\langle {\boldsymbol{\alpha}} \circ {\mathfrak{ F}}^{(j)}_{\mathbf{ c}}, \ \lambda \right\rangle } 
 \ d\nu\left[{\mathbf{ c}}\right]\right|\\ $
$\displaystyle \leq \
\left|\int_{{\mathbf{ D}}_j}
  {\boldsymbol{\phi }} \circ {\mathfrak{ H}}^j \left({\mathbf{ c}}\right) \cdot
{\left\langle {\boldsymbol{\alpha}} \circ {\mathfrak{ F}}^{(j)}_{\mathbf{ c}}, \ \lambda \right\rangle } 
 \ d\nu\left[{\mathbf{ c}}\right]\right|
\ + \
\left|\int_{{\mathcal{ C}}^{\mathbb{M}}\setminus{\mathbf{ D}}_j}
  {\boldsymbol{\phi }} \circ {\mathfrak{ H}}^j \left({\mathbf{ c}}\right) \cdot
{\left\langle {\boldsymbol{\alpha}} \circ {\mathfrak{ F}}^{(j)}_{\mathbf{ c}}, \ \lambda \right\rangle } 
 \ d\nu\left[{\mathbf{ c}}\right]\right|\\ $
$\displaystyle \leq \ 
\int_{{\mathbf{ D}}_j}
 \left|\rule[-0.5em]{0em}{1em} {\boldsymbol{\phi }} \circ {\mathfrak{ H}}^j \left({\mathbf{ c}}\right) \right|\cdot
\left|{\left\langle {\boldsymbol{\alpha}} \circ {\mathfrak{ F}}^{(j)}_{\mathbf{ c}}, \ \lambda \right\rangle }  \right| \ d\nu\left[{\mathbf{ c}}\right]
\ + \
\int_{{\mathcal{ C}}^{\mathbb{M}}\setminus{\mathbf{ D}}_j}
 \left|\rule[-0.5em]{0em}{1em} {\boldsymbol{\phi }} \circ {\mathfrak{ H}}^j \left({\mathbf{ c}}\right) \cdot
{\left\langle {\boldsymbol{\alpha}} \circ {\mathfrak{ F}}^{(j)}_{\mathbf{ c}}, \ \lambda \right\rangle } 
 \right| \ d\nu\left[{\mathbf{ c}}\right] \\ $
$\displaystyle \leq \
\int_{{\mathbf{ D}}_j}
 1 \cdot \frac{\epsilon}{2} \ d\nu\left[{\mathbf{ c}}\right]
\ + \
\int_{{\mathcal{ C}}^{\mathbb{M}}\setminus{\mathbf{ D}}_j} 1 \ d\nu\left[{\mathbf{ c}}\right]
\ \ < \ \ \frac{\epsilon}{2} \ + \ \frac{\epsilon}{2} \ \ = \ \ \epsilon$.
{\tt \hrulefill $\Box$ }\end{list}  \medskip

  Suppose ${\mathcal{ B}}$ is nilpotent, with upper central series (\ref{composition.series}).  If $k\in{\left[ 1...K \right]}$, then ${\mathcal{ Q}}_k = {\mathcal{ Z}}_k/{\mathcal{ Z}}_{k-1}$
is abelian, and the decomposition ${\mathcal{ B}} =
 {\mathcal{ Q}}_1\star\left(
 {\mathcal{ Q}}_2\star\left[\ldots\left({\mathcal{ Q}}_{K-1}\star{\mathcal{ Q}}_K\right)\ldots\right]\right)$
induces a natural identification
${\mathcal{ B}}^{\mathbb{M}}\cong {\mathcal{ Q}}_1^{\mathbb{M}}\times {\mathcal{ Q}}_2^{\mathbb{M}} \times\ldots
\times{\mathcal{ Q}}_K^{\mathbb{M}}$.  If $\lambda_k\in\ensuremath{{\mathcal{ M}}\left[{\mathcal{ Q}}_k^{\mathbb{M}}\right] }$
for all $k$, then
$\lambda_1\otimes\ldots\otimes\lambda_{K}\in\ensuremath{{\mathcal{ M}}\left[{\mathcal{ B}}^{\mathbb{M}}\right] }$.
Let ${\mathcal{ H}}{\mathcal{ M}}\left[{\mathcal{ B}}^{\mathbb{M}}\right]$ denote the convex, weak*-closure in
$\ensuremath{{\mathcal{ M}}\left[{\mathcal{ B}}^{\mathbb{M}}\right] }$ of the set
\[
{\left\{ \lambda_1\otimes\ldots\otimes\lambda_{K} \; ; \; \lambda_k\in{\mathcal{ H}}{\mathcal{ M}}\left[{\mathcal{ Q}}_k^{\mathbb{M}}\right] \ \mathrm{for \ all} \ k \right\} }.
\]

\begin{thm}\label{nilpotent.to.haar}  
  Suppose ${\mathcal{ B}}$ is nilpotent, and ${\mathfrak{ G}}$ has a local map of the
form $\displaystyle{\mathfrak{ g}}({\mathbf{ b}})= b_{{\mathsf{ v}}_1}^{n_1}
 b_{{\mathsf{ v}}_2}^{n_2}\cdots b_{{\mathsf{ v}}_J}^{n_J}$.  For every
${\mathsf{ v}}\in{\mathbb{V}}$, suppose $\displaystyle\ell_{\mathsf{ v}}=\sum_{{\mathsf{ v}}_j={\mathsf{ v}}} n_j$ is
relatively prime to ${{\sf card}\left[{\mathcal{ B}}\right]}$.

  If $\mu\in {\mathcal{ H}}{\mathcal{ M}}\left[{\mathcal{ B}}^{\mathbb{M}}\right]$, then $\displaystyle\mathrm{wk^*}\!\!\lim_{{\mathbb{J}}\ni j{\rightarrow}{\infty}}
{\mathfrak{ G}}^j(\mu) = \eta_{_{{\mathcal{ B}}}}{}$ along a set ${\mathbb{J}}\subset{\mathbb{N}}$ of density one.
Thus, equation (\ref{eqn:cesaro.to.haar}) holds. 
  \end{thm}
\bprf
We'll prove this by induction on $K$, the length of the series
(\ref{composition.series}).  If $K=1$, then ${\mathcal{ B}}$ is abelian;
then ${\mathfrak{ G}}$ is diffusive in density by Example \ref{X:LCA.diffuse},
and the result follows from Theorem 12 of \cite{PivatoYassawi1}.

If $K>1$, then let ${\mathcal{ A}}={\mathcal{ Q}}_1 = Z({\mathcal{ B}})$, 
and ${\mathcal{ C}}={\mathcal{ B}}/{\mathcal{ A}}$.  Thus ${\mathcal{ H}}{\mathcal{ M}}\left[{\mathcal{ B}}^{\mathbb{M}}\right]$ is the
convex weak* closure of ${\mathcal{ S}}={\left\{ \lambda\otimes\nu \; ; \; \lambda\in{\mathcal{ H}}{\mathcal{ M}}\left[{\mathcal{ A}}^{\mathbb{M}}\right] \ \mathrm{and} \ 
\nu\in{\mathcal{ H}}{\mathcal{ M}}\left[{\mathcal{ C}}^{\mathbb{M}}\right] \right\} }$, so it suffices to prove
the theorem for $\mu=\lambda\otimes\nu\in{\mathcal{ S}}$.  Let
${\mathfrak{ G}}={\mathfrak{ F}}\star{\mathfrak{ H}}$.  Then ${\mathfrak{ H}}$ has local map 
 $\displaystyle{\mathfrak{ h}}({\mathbf{ c}})=  c_{{\mathsf{ v}}_1}^{n_1}
  c_{{\mathsf{ v}}_2}^{n_2}\cdots  c_{{\mathsf{ v}}_J}^{n_J}$, and,
by hypothesis, all $\ell_{\mathsf{ v}}$ are all relatively prime to ${{\sf card}\left[{\mathcal{ C}}\right]}$.
But ${\mathcal{ C}}$ has an upper
central series like (\ref{factored.composition.series}) of length
$K-1$, so by induction hypothesis, there is a set
${\mathbb{K}}\subset{\mathbb{N}}$ of density one so that $\displaystyle\mathrm{wk^*}\!\!\lim_{{\mathbb{K}}\ni
k{\rightarrow}{\infty}} {\mathfrak{ H}}^k(\nu) = \eta_{_{{\mathcal{ C}}}}$.

  Since ${\mathcal{ A}}= Z({\mathcal{ B}})$, let
${\mathfrak{ L}}$ be as in {\bf Part 3} of Proposition \ref{MCA.structure.prop}.
As in Example \ref{relCAx0}, $\displaystyle {\mathfrak{ l}}\left({\mathbf{ a}}\right) \ = \
\sum_{{\mathsf{ v}}\in{\mathbb{V}}} \ell_{\mathsf{ v}}\cdot  a_{\mathsf{ v}}$, and, by hypothesis,
$\ell_{\mathsf{ v}}$ are all relatively prime to ${{\sf card}\left[{\mathcal{ A}}\right]}$, so, as in Example
\ref{X:LCA.diffuse}, ${\mathfrak{ L}}$ is ${\mathbb{I}}$-diffusive for some subset
${\mathbb{I}}\subset{\mathbb{N}}$ of density one.  Proposition
\ref{X:aff.diffuse} then implies that ${\mathfrak{ F}}$ is relatively
${\mathbb{I}}$-diffusive.  Let ${\mathbb{J}}={\mathbb{I}}\cap{\mathbb{K}}$, also a set of density one.
Then apply Proposition \ref{P:Reldiffuse.into.Haar} to conclude that
$\displaystyle \mathrm{wk^*}\!\!\lim_{{\mathbb{J}}\ni j{\rightarrow}{\infty}} {\mathfrak{ G}}^j \mu \ = \
\eta_{_{{\mathcal{ A}}}} \otimes \eta_{_{{\mathcal{ C}}}} = \eta_{_{{\mathcal{ B}}}}{}$.
{\tt \hrulefill $\Box$ }\end{list}  \medskip

  \medskip         \refstepcounter{thm}                     \begin{list}{} 			{\setlength{\leftmargin}{1em} 			\setlength{\rightmargin}{0em}}                     \item {\bf Example \thethm:} Recall ${\mathfrak{ G}}:{\mathbf{ Q}}_8^{\mathbb{Z}}\!\longrightarrow{\mathbf{ Q}}_8^{\mathbb{Z}}\!$ \ from
Example (\ref{X:quat.decomp}).  In this case,
${\mathcal{ A}}\cong{{\mathbb{Z}}_{/2}}$, and ${\mathfrak{ L}}$, having 
local map ${\mathfrak{ l}}\left( a_0, a_1, a_2, a_3\right) \ = \
 a_0+ a_1+ a_2+ a_3$, is
diffusive in density, so ${\mathfrak{ F}}$ is relatively diffusive in density.
Meanwhile, ${\mathcal{ C}} = {{\mathbb{Z}}_{/2}}\oplus{{\mathbb{Z}}_{/2}}$ and ${\mathfrak{ H}}$, with
local map ${\mathfrak{ h}}\left( c_0, c_1, c_2, c_3\right) \ = \
 c_0+ c_1+ c_2+ c_3$, is also diffusive in density.
Hence, equation (\ref{eqn:cesaro.to.haar}) holds for any
$\mu\in{\mathcal{ H}}{\mathcal{ M}}\left[{\mathbf{ Q}}_8^{\mathbb{Z}}\right]$.
  	\hrulefill\end{list}

\section{Conclusion}

  Multiplicative cellular automata over a group ${\mathcal{ B}}$ inherit a
natural structural decomposition from ${\mathcal{ B}}$.  Using this
decomposition, we can compute the measurable entropy of MCA, and
show that a broad class of initial measures converge to the Haar
measure in Ces\`aro  average.  However, many questions remain.  For
example, it is unclear how to show relative diffusion when ${\mathcal{ A}}$ is
not central in ${\mathcal{ B}}$.  Indeed, even non-relative diffusion is mysterious
for noncyclic abelian groups \cite{PivatoYassawi2}.  Also,  computation
of relative entropy will be much more complicated in the case of
`variably permutative' relative CA, such as Example (\ref{entrx2});
\ perhaps this requires some `relative' version of Lyapunov exponents
\cite{ShereshevskyLyapunov,Tisseur}.

 Permutative MCA are a considerable generalization of the
linear cellular automata previously studied, but they are still only a
very special class of permutative cellular automata.  The 
asymptotics of measures for general permutative CA \cite{MaassMartinez} is
still poorly understood.
{\footnotesize

\bibliographystyle{plain}
\bibliography{bibliography}
}

\end{document}